\documentclass[11pt, oneside]{article}   	% use "amsart" instead of "article" for AMSLaTeX format
\usepackage{geometry}                		% See geometry.pdf to learn the layout options. There are lots.
\usepackage{graphicx}				% Use pdf, png, jpg, or eps§ with pdflatex; use eps in DVI mode
\usepackage{color}
\usepackage{amssymb}
\usepackage{amsmath}
\usepackage{amsthm}
\usepackage{bm}
\usepackage{aliases}
\usepackage[utf8]{inputenc}
\usepackage[english]{babel}

\newtheorem{theorem}{Theorem}

\title{A global divergence conforming DG method for hyperbolic conservation laws with divergence constraint}
\author{Praveen Chandrashekar\footnote{TIFR Center for Applicable Mathematics, Bangalore, India. Email: \texttt{praveen@tifrbng.res.in}}}

\date{}							% Activate to display a given date or no date

\begin{document}
\maketitle
%\section{}
%\subsection{}
%\linenumbers

\begin{abstract}
We propose a globally divergence conforming discontinuous Galerkin (DG) method on Cartesian meshes for {\em curl-type hyperbolic conservation} laws based on directly evolving the face and cell moments of the Raviart-Thomas approximation polynomials. The face moments are evolved using a 1-D discontinuous Gakerkin method that uses 1-D and multi-dimensional Riemann solvers while the cell moments are evolved using a standard 2-D DG scheme that uses 1-D Riemann solvers. The scheme can be implemented in a local manner without the need to solve a global mass matrix which makes it a truly DG method and hence useful for explicit time stepping schemes for hyperbolic problems. The scheme is also shown to exactly preserve the divergence of the vector field at the discrete level. Numerical results using second and third order schemes for induction equation are presented to demonstrate the stability, accuracy and divergence preservation property of the scheme.
\end{abstract}
\section{Introduction}

Constraint-preserving approximations are important in the numerical simulation of problems in computational electrodynamics (CED) and magnetohydrodynamics (MHD). The time domain Maxwell equations used in CED for the electric and magnetic fields may be written in non-dimensional units as
\[
\df{\E}{t} - \nabla \times \B = -\J, \qquad \df{\B}{t} + \nabla \times \E = 0
\]
with the constraint that
\[
\nabla\cdot\B = 0, \qquad \nabla\cdot\E = \rho, \qquad \df{\rho}{t} + \nabla\cdot\J = 0
\]
where $\rho$ is the electric charge density and $\J$ is the current which can be related to $\E$, $\B$ through Ohm's Law.  In ideal compressible MHD, the magnetic field is given by the Faraday Law or {\em induction} equation
\[
\df{\B}{t} + \nabla \times \E = 0, \qquad \E = -\vel \times \B
\]
with the constraint
\[
\nabla\cdot\B = 0
\]
where $\vel$ is the velocity of the fluid obtained from solving the compressible Euler equations with a Lorentz force that depends on the magnetic field. The above two sets of problems involve {\em hyperbolic conservation laws} with a constraint on the divergence of some vector field and we will concentrate on this type of problems in the present work. The methods we develop in this paper can be applied to the above two class of problems.

There is a large collection of methods developed to solve problems with some divergence constraint spanning Maxwell equations, the simple induction equation and the full MHD equations, both compressible and incompressible. It recognized that satisfying the constraints can have an implication on the accuracy and stability of the schemes~\cite{Brackbill1980426}, \cite{Toth2000605}. Yee~\cite{yee_numerical_1966} proposed a staggered grid central difference scheme for Maxwell equations which preserves a finite difference approximation of the divergence.  This work showed the importance of staggered storage of variables which has been used in other forms by subsequent researchers. A correction method was used for MHD in~\cite{Brackbill1980426} where the magnetic field is first updated by some standard method and then projected to divergence-free space which however requires the solution of a globally coupled problem. In~\cite{1988ApJ...332..659E}, a constrained transport method is developed for MHD which is based on the Yee scheme. Balsara et al.~\cite{BALSARA2001614}, \cite{0067-0049-151-1-149}, \cite{Balsara20095040}, \cite{Balsara2015687} proposed to reconstruct the magnetic field inside the cells in a divergence-free manner, given the information of the normal components on the faces of the cell. The solution on the faces are evolved either with a finite volume or DG method~\cite{Balsara2017104}. These methods require the use of 1-D and 2-D Riemann solvers, where the 2-D Riemann solver involves four states meeting at the vertices of the cells. Recent work has developed methods to solve such 2-D Riemann problems in the context of CED~\cite{Balsara:2016:HRT:2951994.2952032}, \cite{BALSARA2017}, \cite{balsara2017b} and MHD~\cite{BALSARA20101970}, \cite{Balsara:2012:THR:2365370.2365791}, \cite{Balsara:2014:MHR:2580127.2580622}, \cite{Balsara2015269}, \cite{BALSARA201725}.

A DG scheme based on locally divergence-free approximations has been proposed in~\cite{Cockburn:2004:LDD:1008428.1008437, kroner2005, lishu2005} but these are not globally divergence-free since the normal components are not continuous across the cell faces. Central DG schemes for ideal MHD which are globally divergence-free and also free of any Riemann solvers have been proposed in \cite{LI20114828}, \cite{Li:2012:AOE:2109699.2110058} where the magnetic field is approximated by Brezzi-Douglas-Marini (BDM) polynomials~\cite{Brezzi:1991:MHF:108342} on staggered Cartesian meshes. A stability analysis of the first order central DG scheme when applied to induction equation has been performed in~\cite{refId0}.

In Lagrange multiplier methods, an artificial pressure is introduced in the induction equation and the divergence-free condition is satisfied in a weak sense, see e.g. \cite{schotzau_mixed_2004}. There is also the class of hyperbolic divergence cleaning methods~\cite{Dedner2002645} where the divergence errors are damped by adding an artificial pressure like term. A summation-by-parts finite difference scheme for induction equation have been developed in~\cite{Koley2009}, \cite{doi:10.1093/imanum/drq030} where a term proportional to the divergence is added to the equations to obtain a stable scheme. A stable upwind finite difference scheme based on the symmetrized version of the equations in the non-conservative form is constructed in \cite{fuchs2009}. For MHD, methods have been developed based on Godunov's symmetrized version of the equations where stability is achieved using Riemann solvers~\cite{powell94} or using entropy stability ideas~\cite{doi:10.1137/15M1013626}, \cite{Winters:2016:AEC:2846979.2847024}. But all these approaches are non-conservative since they modify the PDE in a non-conservative manner.

There are also a wide variety of {\em Galerkin} finite element methods developed for curl-type equations like Maxwell, induction and MHD, see e.g., \cite{hiptmair_2002}, \cite{Hu2017}. When applied to MHD problems as in~\cite{Hu2017}, the magnetic field is approximated by a $\hdiv$ conforming space while the electric field is approximated by a $\hcurl$ conforming space. The first type of spaces have continuous normal components while the second type have continuous tangential components across the cell faces. Divergence-free bases have been developed in~\cite{Cai2013}, \cite{cai2017} that allow to locally correct any divergence error. In all of these methods, the approximating spaces require some continuity across the cell faces leading to global matrices which require efficient matrix solution techniques for application to large scale problems arising in real world situations. If we are interested in purely hyperbolic problems like Maxwell equations or ideal MHD where explicit time stepping schemes are used, we will have to solve a global mass matrix in every time step which increases the computational expense. One of the goals of this work is to construct a divergence conforming method that has local mass matrices which can be easily inverted on each cell or face.

Coming back to the present work, as a prototypical model, we will consider the curl-type equation of the form
\begin{equation}
\df{\B}{t} + \nabla \times \E = -\M
\label{eq:ind}
\end{equation}
whose solutions satisfy
\begin{equation}
\df{}{t}(\nabla\cdot\B) + \nabla\cdot \M = 0
\label{eq:diveq}
\end{equation}
If $\M=0$ and $\nabla\cdot\B = 0$ at the initial time, then we have $\nabla\cdot\B = 0$ at future times also. This property does not depend on the particular form of $\E$ but as a concrete example we will take $\E = - \vel \times \B$ as in the induction equation with $\vel$ being a given velocity field. We will also restrict the examples to two dimensional case where the equations are of the form
\[
\df{\Bx}{t} + \df{E}{y} = -M_x, \qquad \df{\By}{t} - \df{E}{x} = -M_y
\]
with $E = \vy \Bx - \vx \By$.

The goal of this paper is to present a divergence constraint preserving DG scheme for the curl type equation of the form~(\ref{eq:ind}) that makes use of the standard Raviart-Thomas elements~\cite{Raviart1977} on Cartesian meshes. The Raviart-Thomas polynomials are defined in terms of certain face and cell moments. In the present work, the face moments are evolved using a DG scheme proposed in Balsara \& Kappeli~\cite{Balsara2017104} which makes use of 1-D and 2-D Riemann solvers. However, unlike those authors, the novelty of our formulation is that it is based from the ground-up on the Raviart-Thomas elements. As a result, our scheme also includes internal nodes within each element in addition to the facial nodes proposed in~\cite{Balsara2017104} but there is no need to perform a divergence-free reconstruction step. The internal nodal values are evolved according to a conventional DG scheme which makes of 1-D Riemann solver to obtain the numerical fluxes. In 2-D, the algorithm involves a one dimensional DG scheme built on the faces of the mesh and a two dimensional DG scheme built in the interior of the cells. The idea of evolving moments of BDM polynomials has been used in the central DG scheme~\cite{LI20114828}, \cite{Li:2012:AOE:2109699.2110058}  together with a reconstruction step and the use of staggered grids, but in our approach there is no need to perform a recontruction step and all the degrees of freedom are directly evolved by the DG scheme. By having certain compatibility in the internal and facial DG schemes, we obtain a scheme that overall preserves the divergence of the solution and also ensures that normal component of the vector field is continuous across the cell faces. This synthesis yields a conceptually pleasing time-explicit DG scheme which avoids staggered meshes and brings together methods for CED and MHD. An important feature of the present scheme is that {\em there is no need to invert a global mass matrix} which makes explicit time stepping to be very efficient. While only second and third order results are shown here, the extension to much higher order is very easy since we can easily construct high order Raviart-Thomas polynomials on Cartesian meshes.  Since it does not require staggered grids, the methodology presented here can be used for other curl-type equations, including the Maxwell equations, to develop constraint preserving schemes on unstructured quadrilateral/hexahedral grids, isoparametric elements and also on adaptively refined grids of quadtree/octree type. An unstaggered DG scheme using BDM polynomials as approximation space has been proposed in~\cite{Fu2018} for induction equation and the full MHD system which is able to preserve the divergence of the magnetic field. This scheme makes use of multi-dimensional Riemann solvers at the cell vertices and evolves the normal components on the faces by a DG scheme, and the present scheme is very similar to~\cite{Fu2018}. A Fourier stability analysis for induction equation in~\cite{Fu2018} shows the importance of correctly approximating the multi-dimensional Riemann solution and the theoretical developments there will be useful in extending our own scheme to the MHD system.

The rest of the paper is organized as follows. The schemes proposed here make use of some non-standard approximation spaces which is not so well known among scientific and engineering community. Hence we have strived to provide some elementary introduction and derivations to explain this important topic to a wider audience. We start in section~(\ref{sec:app}) by explaining the process of approximating vector fields whose divergence has to be bounded in terms of Raviart-Thomas polynomials. The construction of the approximation in terms of moments is explained and we show how the approximation automatically satisfies the divergence-free condition. In section~(\ref{sec:dg}), the DG scheme is proposed to evolve the moments and its ability to exactly preserve the divergence is shown. Then we detail the numerical fluxes and boundary conditions for the induction equation. Finally, in section~(\ref{sec:res}), we show through many numerical tests that the proposed schemes have to optimal accuracy in approximating both divergence-free and divergent solutions.
%-----------------------------------------------------------------------------------------------
\section{Approximation of vector fields}
\label{sec:app}
\begin{figure}
\begin{center}
\includegraphics[width=0.8\textwidth]{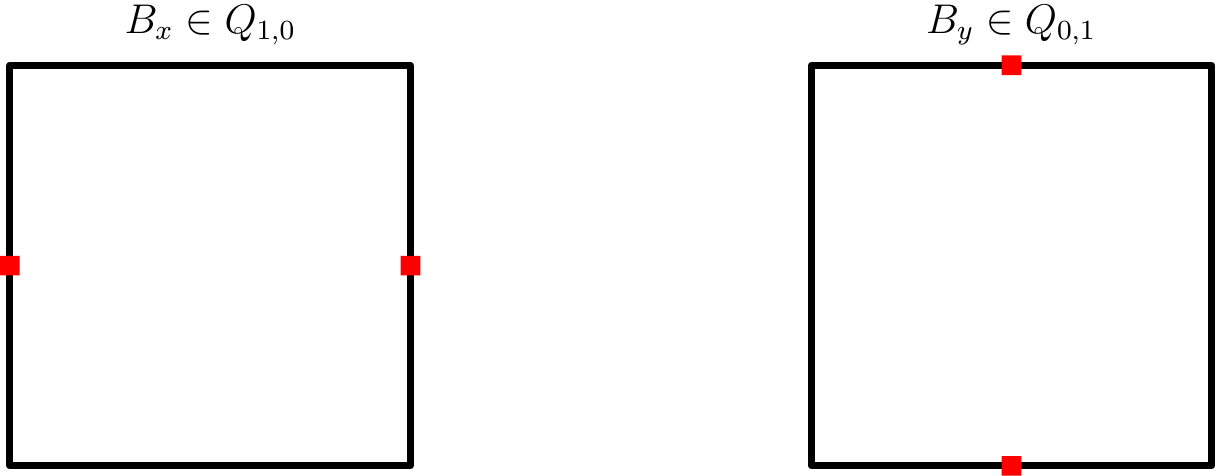}
\caption{Location of dofs of Raviart-Thomas polynomial for $k=0$}
\label{fig:dof0}
\end{center}
\end{figure}

\begin{figure}
\begin{center}
\includegraphics[width=0.8\textwidth]{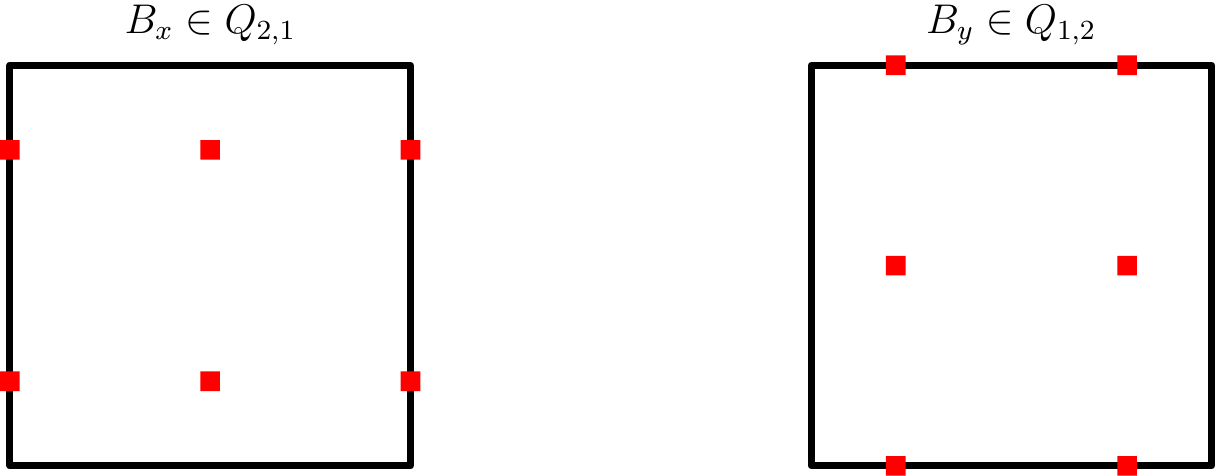}
\caption{Location of dofs of Raviart-Thomas polynomial for $k=1$}
\label{fig:dof1}
\end{center}
\end{figure}

\begin{figure}
\begin{center}
\includegraphics[width=0.8\textwidth]{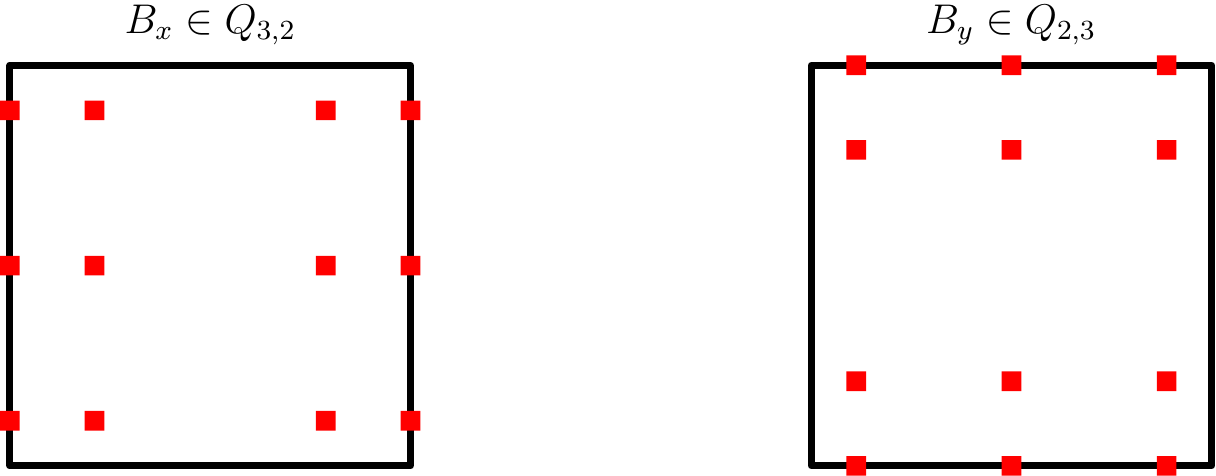}
\caption{Location of dofs of Raviart-Thomas polynomial for $k=2$}
\label{fig:dof2}
\end{center}
\end{figure}

When dealing with problems where the vector field $\B$ must be divergence-free, it is natural to look for solutions in the space $\hdiv$ which is defined as
\[
\hdiv = \{ \B \in \LtOmega : \Div(\B) \in \ltOmega\}
\]
i.e., these functions have bounded $L^2$ norm and the divergence also has bounded $L^2$ norm. To approximate functions in $\hdiv$ on a mesh $\mesh$ with piecewise polynomials as done in the finite element method, we need the following compatibility condition.
\begin{theorem}[See \cite{Quarteroni:2008:NAP:1502402}, Proposition 3.2.2]
Let $\Bh : \Omega \to \re^d$ be such that
\begin{enumerate}
\item $\Bh|_K \in \bm{H}^1(\Omega)$ for all $K \in \mesh$
\item for each common face $F = K_1 \cap K_2$, $K_1, K_2 \in \mesh$, the trace of normal component $\un \cdot \Bh|_{K_1}$ and $\un \cdot \Bh|_{K_2}$ is the same.
\end{enumerate}
Then $\Bh \in \hdiv$. Conversely, if $\Bh \in \hdiv$ and (1) holds, then (2) is also satisfied.
\end{theorem}
The functions in $\hdiv$ can be approximated on a Cartesian mesh by the Raviart-Thomas space of piecewise polynomial functions as follows. Define the one dimensional polynomials $P_k(x)$, $P_k(y)$ of degree at most $k$ with respect to the variables $x$, $y$ respectively. Let $Q_{r,s}(x,y)$ denote the tensor product polynomials of degree $r$ in the variable $x$ and degree $s$ in the variable $y$, i.e.,
\[
Q_{r,s}(x,y) = \vspan\{ x^i y^j, \ 0 \le i \le r, \ 0 \le j \le s \}
\]
For $k \ge 0$, the Raviart-Thomas space of vector functions is defined as
\[
\rt_k = Q_{k+1,k} \times Q_{k,k+1}
\]
The dimension of this space is $2(k+1)(k+2)$. For any $\Bh \in \rt_k$, we have $\Div(\Bh) \in Q_{k,k}(x,y)$. We consider a cell centered at the origin and of size $\dx, \dy$. The restriction of $\Bh = (\Bxh, \Byh)$ to a face is a polynomial of degree $k$, i.e.,
\[
\Bxh(\pm\dx/2,y) \in P_k(y), \qquad \Byh(x,\pm\dy/2) \in P_k(x)
\]
For doing the numerical computations, it is useful to map each cell $C$ to a reference cell and choose certain nodes that can be used to define Lagrange polynomials.  Let $\{\xi_i, 0 \le i \le k+1\}$  and $\{\hxi_i, 0 \le i \le k\}$ be two sets of distinct nodes in the reference interval $[0,1]$ with the constraint that $\xi_0 = 0$ and $\xi_{k+1}=1$. Let $\phi_i$ and $\hphi_i$ be the corresponding one dimensional Lagrange polynomials. Then the magnetic field is given by
\begin{equation}
\Bxh(\xi,\eta) = \sum_{i=0}^{k+1} \sum_{j=0}^k (\Bx)_{ij} \phi_i(\xi) \hphi_j(\eta), \qquad \Byh(\xi,\eta) = \sum_{i=0}^{k} \sum_{j=0}^{k+1} (\By)_{ij} \hphi_i(\xi) \phi_j(\eta)
\label{eq:Blag}
\end{equation}
Our choice of nodes ensures that the normal component of the magnetic field is continuous on the cell faces. There is no unique way to choose the nodes and the particular choice we use is what is implemented in the \texttt{deal.II} library~\cite{BangerthHartmannKanschat2007}. We first describe the nodes used for the $B_x$ component. On the left and right faces of the cell, we have $k+1$ Gauss-Legendre nodes, while on the interior, we have tensor product of $k \times (k+1)$ Gauss-Legendre nodes. For the $B_y$ component, we have $(k+1)$ Gauss-Legendre nodes on the bottom and top faces of the cell, and a tensor product of $(k+1) \times k$ Gauss-Legendre nodes in the interior. Figures~(\ref{fig:dof0}), (\ref{fig:dof1})  and (\ref{fig:dof2}) show the location of the nodes for the case of $k=0$, $k=1$ and $k=2$ respectively. On Cartesian meshes, the optimal error estimates for approximating vector fields with $\rt_k$ are of the form~\cite{Brezzi:1991:MHF:108342}, \cite{doi:10.1137/S0036142903431924}
\begin{eqnarray}
\label{eq:app1}
\norm{\B - \Bh}_{\LtOmega} &\le& C h^{k+1} |\B|_{\HOmega} \\
\label{eq:app2}
\norm{\Div(\B) - \Div(\Bh)}_{\LtOmega} &\le& C h^{k+1} |\Div(\B)|_{\HOmega}
\end{eqnarray}
The second error estimate implies that $\Div(\Bh) \equiv 0$ if $\Div(\B) \equiv 0$ but we will show the divergence-free property of the approximation more explicitly in a later section.

%-----------------------------------------------------------------------------------------------
\subsection{Construction of $\B_h$}
\begin{figure}
\begin{center}
\includegraphics[width=0.4\textwidth]{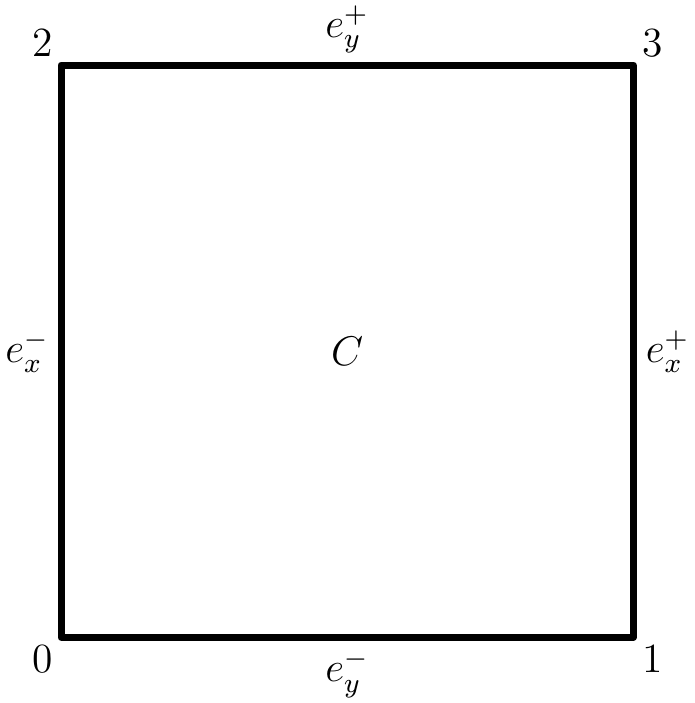}
\caption{A cell $C$ and its face nomenclature}
\label{fig:cell}
\end{center}
\end{figure}
To determine a function $\Bh \in \rt_k$ on each cell $C$, we need $2(k+1)(k+2)$ pieces of information which is the dimension of the $\rt_k$ space. These are taken to be certain moments on the faces and interior of the cell. The face moments are given by
\[
\int_{\expm} \Bxh \phi \ud y \qquad \forall \phi \in P_k(y)
\]
and
\[
\int_{\eypm} \Byh \phi \ud x \qquad \forall \phi \in P_k(x)
\]
where $\expm$ are the two vertical sides of cell $C$ and $\eypm$ are the two horizontal sides of cell $C$ as shown in figure~(\ref{fig:cell}). The cell moments are given by
\[
\int_C \Bxh \psi \ud x \ud y \qquad \forall \psi \in \partial_x Q_{k,k}(x,y) := Q_{k-1,k}(x,y)
\]
and
\[
\int_C \Byh \psi \ud x \ud y \qquad \forall \psi \in \partial_y Q_{k,k}(x,y) := Q_{k,k-1}(x,y)
\]
Note that $ \dim P_k(x) = \dim P_k(y) = k+1$ and $\dim \partial_x Q_{k,k}(x,y) = \dim \partial_y Q_{k,k}(x,y) = k(k+1)$ so that we have in total $4(k+1)+2k(k+1) = 2(k+1)(k+2)$ pieces of information which is enough to determine $\Bh \in \rt_k$. The moments on the faces $\expm$ uniquely determine the restriction of $\Bxh$ on those faces, and similarly the moments on $\eypm$ uniquely determine the restriction of $\Byh$ on the corresponding faces. This ensures continuity of the normal component of $\Bh$ on all the faces. Also note that the moment equations for $\Bxh$, $\Byh$ are decoupled and can be solved independently of one another. Given all the moments for a cell, we can uniquely determine the function $\B_h$ as shown in the following theorem. If all the moments are zero, then we show that the $\B_h$ is identically zero which implies that the matrix arising from the moments is invertible. This is a standard result which we state here in a simple setting of Cartesian meshes and for a proof of this result on general meshes, see~(\cite{Brezzi:1991:MHF:108342}, Proposition 3.3 and 3.4) and~\cite{Nedelec1980}.
%-----------------------------------------------------------------------------------------------
\begin{theorem}
If all the moments are zero for any cell $C$, then $\Bh \equiv 0$ inside that cell.
\end{theorem}
\noindent
\underline{Proof}: The face moments being zero implies that 
\[
\Bxh \equiv 0 \quad \textrm{on} \quad \expm \qquad\textrm{and}\qquad \Byh \equiv 0 \quad \textrm{on} \quad \eypm
\]
Now take $\psi = \partial_x \phi$ for some $\phi \in Q_{k,k}$ in the cell moment equation of $\Bxh$ and perform an integration by parts
\[
-\int_C \df{\Bxh}{x} \phi \ud x \ud y - \int_{\exl} \Bxh \phi \ud y + \int_{\exr} \Bxh \phi \ud y = 0
\]
and hence
\[
\int_C \df{\Bxh}{x} \phi \ud x \ud y = 0 \qquad \forall \phi \in Q_{k,k}
\]
Since $\df{\Bxh}{x} \in Q_{k,k}$, this implies that $\df{\Bxh}{x} \equiv 0$ and hence $\Bxh \equiv 0$. Similarly, we conclude that $\Byh \equiv 0$. \qed
%-----------------------------------------------------------------------------------------------
\begin{theorem}
Let $\Bh \in \rt_k$ satisfy the moments
\begin{eqnarray}
\label{eq:mom1}
\int_{\expm} \Bxh \phi \ud y &=& \int_{\expm} \Bx \phi \ud y \qquad \forall \phi \in P_k(y) \\
\label{eq:mom2}
\int_{\eypm} \Byh \phi \ud x &=& \int_{\eypm} \By \phi \ud x \qquad \forall \phi \in P_k(x) \\
\label{eq:mom3}
\int_C \Bxh \psi \ud x \ud y &=& \int_C \Bx \psi \ud x \ud y \qquad \forall \psi \in \partial_x Q_{k,k}(x,y) \\
\label{eq:mom4}
\int_C \Byh \psi \ud x \ud y &=& \int_C \By \psi \ud x \ud y \qquad \forall \psi \in \partial_y Q_{k,k}(x,y)
\end{eqnarray}
for a given vector field $\B \in \hdiv$. If $\Div (\B) \equiv 0$ then $\Div (\Bh) \equiv 0$.
\end{theorem}
\noindent
\underline{Proof}: This is also a standard result and we refer the reader to the textbook~\cite{Brezzi:1991:MHF:108342} for a general proof and other references. We choose $\psi = \partial_x \phi$ and $\psi=\partial_y\phi$ for some $\phi \in Q_{k,k}(x,y)$ respectively in the two cell moment equations (\ref{eq:mom3}), (\ref{eq:mom4}). Adding these two equations together, we get
\[
\int_C (\Bxh \partial_x\phi + \Byh \partial_x\phi) \ud x \ud y = \int_C (\Bx \partial_x\phi + \By \partial_y \phi) \ud x \ud y
\]
Performing integration by parts on both sides
\[
- \int_C (\partial_x \Bxh + \partial_y \Byh ) \phi \ud x \ud y + \int_{\partial C} \phi \Bh \cdot \un \ud s = -\int_C (\partial_x \Bx + \partial_y \By) \phi \ud x \ud y + \int_{\partial C} \phi \B \cdot \un \ud s
\]
where $\un$ is the unit normal vector on the boundary of the cell. Note that $\phi$ restricted to $\partial C$ is a one dimensional polynomial of degree $k$ and the face moments of $\Bh$ and $\B$ agree with one another by equations~(\ref{eq:mom1}), (\ref{eq:mom2}). Hence we get
\[
\int_C (\partial_x \Bxh + \partial_y \Byh ) \phi \ud x \ud y  = \int_C (\partial_x \Bx + \partial_y \By) \phi \ud x \ud y \qquad \forall \phi \in Q_{k,k}(x,y)
\]
If $\Div(\B) \equiv 0$, then
\[
\int_C \Div(\Bh) \phi \ud x \ud y  = 0 \qquad \forall \phi \in Q_{k,k}(x,y)
\]
Since $\Div(\Bh) \in Q_{k,k}(x,y)$ this implies that $\Div(\Bh) \equiv 0$ everywhere inside the cell $C$. \qed
%-----------------------------------------------------------------------------------------------
\paragraph{Remark} The moments are with respect to some non-standard test function spaces whose choice is now well motivated by the above theorem. Since the test functions for the cell moments can be written as $\partial_x \phi$ and $\partial_y \phi$ for some $\phi \in Q_{k,k}$, we can obtain an equation for the divergence by doing an integration by parts, which helps us to control the divergence of the approximation.
%-----------------------------------------------------------------------------------------------
\paragraph{Remark} The proof makes use of integration by parts for which the quadrature must be exact. The integrals involving $\Bh$ can be evaluated exactly using Gauss quadrature of sufficient accuracy. This is not the case for the integrals involving $\B$ since it can be an arbitrary nonlinear function. When $\Div(\B)=0$, we have $\B = (\partial_y \Phi, -\partial_x \Phi)$ for some smooth function $\Phi$. We can approximate $\Phi$ by $\Phi_h \in Q_{k+1,k+1}$ and compute the projections using $(\partial_y \Phi_h, -\partial_x \Phi_h)$ in which case the integrations can be performed exactly. This procedure is used in all the test cases with zero divergence which then ensures that the initial condition $\Bh(0)$ has zero divergence everywhere.
%-----------------------------------------------------------------------------------------------
\section{DG scheme for the induction equation}
\label{sec:dg}
We propose a mixed scheme which evolves the face and cell moments that are used to define the projection on the Raviart-Thomas space. Since we have shown that this projection gives divergence-free approximations, the evolution of the same moments will be able to preserve the divergence of the magnetic field at future times also. The face moments are evolved using the scheme proposed in Balsara \& Kappeli~\cite{Balsara2017104}. The variation of $B_x$ on the vertical faces $\expm$ is given by a one dimensional PDE in the $y$ direction, and similarly for the $B_y$ component. Hence we can discretize the one dimensional PDEs on the faces by applying a DG method. Multiplying by test functions used to define the face moments and performing an integration by parts, the DG scheme on the faces is given by
\begin{equation}
\int_{\expm} \df{\Bxh}{t} \phi \ud y - \int_\expm \hE \df{\phi}{y} \ud y + [\tE \phi]_\expm = -\int_{\expm} \hM_x \phi \ud y \qquad \forall \phi \in P_k(y)
\label{eq:dg1}
\end{equation}
\begin{equation}
\int_{\eypm} \df{\Byh}{t} \phi \ud x + \int_\eypm \hE \df{\phi}{x} \ud x - [\tE \phi]_\eypm = -\int_{\eypm} \hM_y \phi \ud x \qquad \forall \phi \in P_k(x)
\label{eq:dg2}
\end{equation}
where $\hE$ is a numerical flux from a 1-D Riemann solver which is required on the faces of the cells, while $\tE$ is a numerical flux obtained from a multi-D Riemann solver and is needed at the vertices of the cells. Note that we may have to approximate $M_x, M_y$ on the faces by some numerical scheme if the source term depends on the solution, since the tangential components of $\Bh$ could be discontinuous across the cell faces, and this numerical approximation is denoted as $\hM_x,\hM_y$, respectively. This type of situation occurs for Maxwell equations as discussed in the Introduction. The quantities $[\tE\phi]_\expm$, $[\tE\phi]_{\eypm}$ are difference operators on the faces, and with respect to the vertex numbering in figure~(\ref{fig:cell}), are defined as follows:
\[
[\tE \phi]_\exl = (\tE\phi)_2 - (\tE\phi)_0, \qquad [\tE \phi]_\exr = (\tE\phi)_3 - (\tE\phi)_1
\]
\[
[\tE \phi]_\eyb = (\tE\phi)_1 - (\tE\phi)_0, \qquad [\tE \phi]_\eyt = (\tE\phi)_3 - (\tE\phi)_2
\]
The cells moments are evolved by the following standard DG scheme
\begin{equation}
\int_C \df{\Bxh}{t} \psi \ud x \ud y - \int_C E \df{\psi}{y} \ud x \ud y + \int_{\partial C} \hE \psi n_y \ud s = -\int_C M_x \psi \ud x \ud y \qquad \forall \psi \in \partial_x Q_{k,k}(x,y)
\label{eq:dg3}
\end{equation}
\begin{equation}
\int_C \df{\Byh}{t} \psi \ud x \ud y + \int_C E \df{\psi}{x} \ud x \ud y - \int_{\partial C} \hE \psi n_x \ud s = -\int_C M_y \psi \ud x \ud y \qquad \forall \psi \in \partial_y Q_{k,k}(x,y)
\label{eq:dg4}
\end{equation}
Note that the same 1-D numerical flux $\hE$ is used in both the face and cell moment equations whereas the vertex numerical flux $\tE$ is needed only in the face moment equations. 

%-----------------------------------------------------------------------------------------------
\begin{theorem}
If $\M=0$, the DG scheme~(\ref{eq:dg1})-(\ref{eq:dg4}) preserves the divergence of the magnetic field. If $\M \ne 0$, then the divergence evolves consistently with equation~(\ref{eq:diveq}) in the sense that the numerical solution satisfies
\begin{equation}
\int_C \phi \df{}{t}\Div(\Bh) \ud x \ud y - \int_C \M \cdot \nabla\phi \ud x \ud y + \int_{\partial C} \phi \bm{\hat{M}} \cdot \un \ud s = 0, \qquad \forall \phi \in Q_{k,k}
\label{eq:divdg}
\end{equation}
\end{theorem}
\noindent
\underline{Proof}: For any $\phi \in Q_{k,k}(x,y)$ take test functions $\psi = \partial_x \phi$ and $\psi = \partial_y \phi$ in the two cell moment equations~(\ref{eq:dg3}),~(\ref{eq:dg4}) respectively and add them together to obtain
\[
\begin{aligned}
\int_C \left[ \df{\Bxh}{t} \partial_x \phi + \df{\Byh}{t} \partial_y \phi \right] \ud x \ud y + & \int_C \M \cdot \nabla\phi \ud x \ud y
-  \int_{\partial C}\hE (n_x \partial_y \phi - n_y \partial_x \phi) \ud s = 0
\end{aligned}
\]
Note that two of the cell integrals cancel since $\partial_x \partial_y \phi = \partial_y \partial_x \phi$. Performing an integration by parts in the first term, we obtain
\begin{equation}
\begin{aligned}
- \int_C \phi \df{}{t}(\partial_x \Bxh +  \partial_y \Byh) \ud x \ud y + & \int_{\partial C} \phi \df{}{t}(\Bh \cdot \un) \ud s \\
+ & \int_C \M \cdot \nabla\phi \ud x \ud y - \int_{\partial C}\hE (n_x \partial_y \phi - n_y \partial_x \phi) \ud s = 0
\end{aligned}
\label{eq:pe1}
\end{equation}
Now, let us concentrate on the second and last terms which can be re-arranged as follows
\begin{eqnarray*}
&& \int_{\partial C} \phi \df{}{t}(\Bh \cdot \un) \ud s - \int_{\partial C}\hE (n_x \partial_y \phi - n_y \partial_x \phi) \ud s \\
&=& \int_{\exr} \phi \df{\Bxh}{t} \ud y - \int_{\exl} \phi \df{\Bxh}{t} \ud y + \int_{\eyt} \phi \df{\Byh}{t} \ud x - \int_{\eyb} \phi \df{\Byh}{t} \ud x \\
&& - \int_{\exr} \hE \partial_y \phi \ud y + \int_{\exl} \hE \partial_y \phi \ud y + \int_{\eyt} \hE \partial_x \phi \ud x - \int_{\eyb} \hE \partial_x \phi \ud x
\end{eqnarray*}
The restriction of $\phi$ on each face is a one dimensional polynomial of degree $k$ and we can use the face moment equations~(\ref{eq:dg1}),~(\ref{eq:dg2}) to obtain
\begin{eqnarray}
\nonumber
&& \int_{\partial C} \phi \df{}{t}(\Bh \cdot \un) \ud s - \int_{\partial C}\hE (n_x \partial_y \phi - n_y \partial_x \phi) \ud s \\
\nonumber
&=& -\int_{\exr} \hM_x \phi \ud y + \int_{\exl} \hM_x \phi \ud y - \int_{\eyt} \hM_y \phi \ud x + \int_{\eyb} \hM_y \phi \ud x\\
\nonumber
&& - [\tE \phi]_{\exr} + [\tE \phi]_{\exl} + [\tE \phi]_{\eyt} - [\tE \phi]_{\eyb} \\
&=& -\int_{\partial C} \phi \bm{\hat{M}} \cdot \un \ud s 
\label{eq:pe2}
\end{eqnarray}
since
\begin{eqnarray*}
&&- [\tE \phi]_{\exr} + [\tE \phi]_{\exl} + [\tE \phi]_{\eyt} - [\tE \phi]_{\eyb} \\
&=& -(\tE\phi)_3 + (\tE\phi)_1 + (\tE\phi)_2 - (\tE\phi)_0 + (\tE\phi)_3 - (\tE\phi)_2 - (\tE\phi)_1 + (\tE\phi)_0 \\
&=& 0
\end{eqnarray*}
Combining equations~(\ref{eq:pe1}) and (\ref{eq:pe2}), we obtain equation~(\ref{eq:divdg}). In the case of $\M=0$, any consistent numerical approximation would lead to $\hat{\M}=0$, and then we obtain
\[
\int_C \phi \df{}{t}(\Div(\Bh)) \ud x \ud y = 0 \qquad \forall \phi \in Q_{k,k}(x,y)
\]
Since $\Div(\Bh) \in Q_{k,k}(x,y)$, we conclude that the divergence is preserved by the numerical scheme.  \qed
%-----------------------------------------------------------------------------------------------
\paragraph{Remark} The above proof required integration by parts in the terms involving the time derivative which is usually called the mass matrix. The other cell integral in the DG scheme can be computed using any quadrature rule of sufficient order and need not be exact. All the face integrals which involve the numerical flux $\hE$ appearing in the face moment and cell moment evolution equations must be computed with the same rule and it is not necessary to be exact for the above proof to hold. However, from an accuracy point of view, these quadratures must be of a sufficiently high order to obtain optimal error estimates. In practice we find that using a $(k+1)$-point Gauss-Legendre quadrature for face integrals and a tensor product rule of the same points for the cell integrals leads to optimal convergence rates. 
%-----------------------------------------------------------------------------------------------
\paragraph{Remark} The preservation of divergence does not rely on the specific form of the fluxes $\tE$, $\hE$ but only on the fact that we have a unique flux $\tE$ at all the vertices, and that we use the same 1-D numerical flux $\hE$ in both the face and cell moment equations.
%-----------------------------------------------------------------------------------------------
\paragraph{Remark}
In the case of Maxwell equations, the electric field has a curl form just like the induction equation but can also have a source term related to the electric current. In our notation, this would correspond to the case when the source term $\M \ne 0$. As discussed in the introduction, we would like to compute the solution and also look at its divergence since it gives information about the charge density in space. From the proof of the previous theorem, we have seen that the divergence satisfies equation~(\ref{eq:divdg}) that looks like a standard DG scheme for equation~(\ref{eq:diveq}). The divergence is a tensor product polynomial of degree $k$, i.e., $\Div(\Bh) \in Q_{k,k}$, and we can expect $\Div(\Bh)$ to be accurate to $O(h^{k+1})$. This is indeed borne out in our numerical tests which shows that we can compute the charge density to the same order of accuracy as the solution without any extra effort.

%-----------------------------------------------------------------------------------------------
\subsection{Numerical fluxes and Electric Fields}
In order to complete the description of the DG scheme, we have make a choice of the two types of numerical fluxes needed in the scheme. To specify the numerical fluxes, we have to identify the characteristic curves in the PDE. Using the zero divergence condition, we can rewrite the induction equation in the following way
\[
\df{\Bx}{t} + \vel \cdot \nabla \Bx + \Bx \df{\vy}{y} - \By \df{\vx}{y} = 0, \qquad \df{\By}{t} + \vel \cdot \nabla \By + \By \df{\vx}{x} - \Bx \df{\vy}{x} = 0
\]
There is only one set of characteristic curves and they are the integral curves of $\vel$ and the velocity field is assumed to be given as a function of space and time coordinates. Following the upwind principle that information propagates along characteristics, the 1-D numerical flux is given by
\[
\hE = \begin{cases}
E_L & \textrm{if } \vel \cdot \un > 0 \\
E_R & \textrm{otherwise}
\end{cases}
\]
where the subscripts L and R denote the left and right states with the normal vector $\un$ pointing from L to R. For example, across the face $\expm$, the flux is given by
\[
\hE = \begin{cases}
\vy \Bx - \vx \By^L & \textrm{if } \vx > 0 \\
\vy \Bx - \vx \By^R & \textrm{otherwise}
\end{cases}
\]
\begin{figure}
\begin{center}
\includegraphics[width=0.45\textwidth]{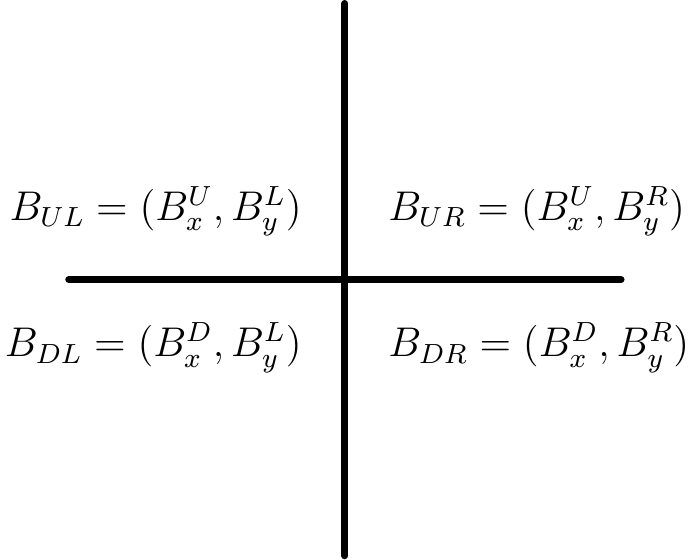}
\caption{Four states meeting at a vertex that define a 2-D Riemann problem}
\label{fig:corner}
\end{center}
\end{figure}
At a vertex, we have four states meeting which is illustrated in figure~(\ref{fig:corner}). Note that $\Bx$ is continuous across the vertical faces and $\By$ is continuous across the horizontal faces. The upwinded electric field at the vertices of the two-dimensional mesh is given by
\[
\tE = \begin{cases}
E_{DL} & \textrm{if } \vx > 0, \ \vy > 0 \\
E_{UL} & \textrm{if } \vx > 0, \ \vy < 0 \\
E_{DR} & \textrm{if } \vx < 0, \ \vy > 0 \\
E_{UR} & \textrm{if } \vx < 0, \ \vy < 0
\end{cases}
\]
which can be written in compact form as
\[
\tE =  \frac{\vy}{2}(\Bx^{U} + \Bx^{D}) - \frac{\vx}{2}(\By^{L} + \By^{R}) 
 - \frac{|\vy|}{2}\left( \Bx^{U} - \Bx^{D}\right) + \frac{|\vx|}{2}\left( \By^{R} -\By^{L} \right)
\]
An equivalent expression is given by~\cite{Balsara2017104}
\begin{equation}
\begin{aligned}
\tE = & \frac{\vy}{4}(\Bx^{UL} + \Bx^{UR} + \Bx^{DL} + \Bx^{DR}) - \frac{\vx}{4}(\By^{UL} + \By^{UR} + \By^{DL} + \By^{DR}) \\
& - \frac{|\vy|}{2}\left( \frac{\Bx^{UL} + \Bx^{UR}}{2} - \frac{\Bx^{DL} + \Bx^{DR}}{2} \right) + \frac{|\vx|}{2}\left( \frac{\By^{UR} + \By^{DR}}{2} - \frac{\By^{UL} + \By^{DL}}{2} \right)
\end{aligned}
\label{eq:nfbk}
\end{equation}
with the understanding that $\Bx^{DL}=\Bx^{DR}$, etc.  We refer the reader to~\cite{Fu2018} for a stability analysis of the first order scheme with the above numerical fluxes. Of course, in the system case like full MHD, the expressions are not so simple and we point the reader to the work in~\cite{BALSARA20101970}, \cite{Balsara:2012:THR:2365370.2365791}, \cite{Balsara:2014:MHR:2580127.2580622}, \cite{Balsara2015269}.
%-----------------------------------------------------------------------------------------------
\subsection{Boundary condition}
\begin{figure}
\begin{center}
\includegraphics[width=0.8\textwidth]{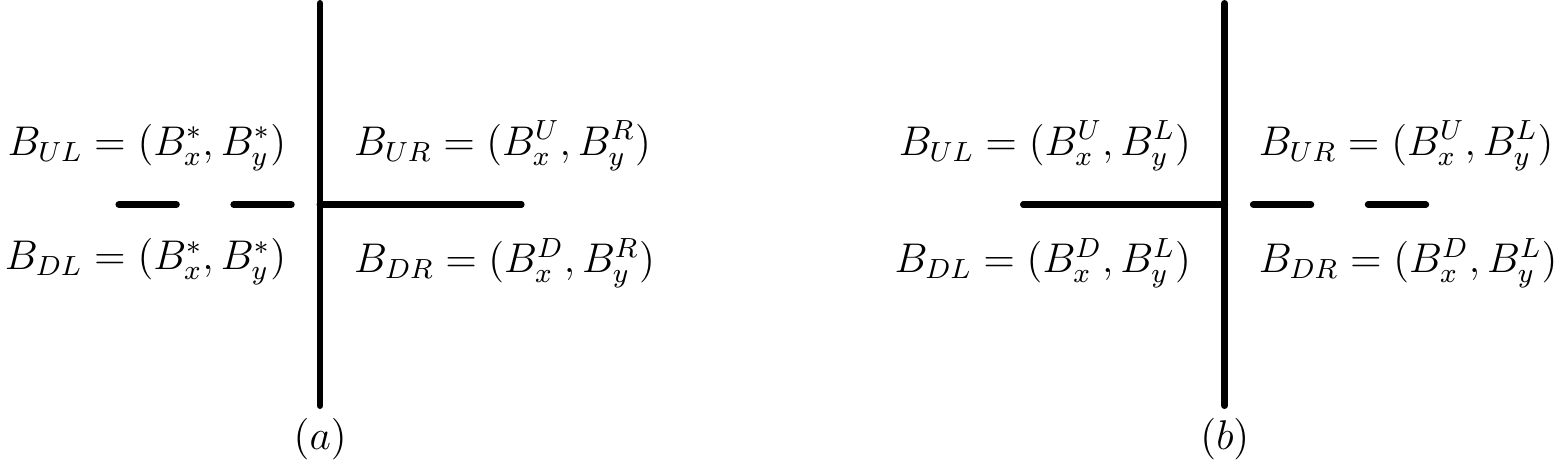}
\caption{Vertex states at boundary: (a) inflow vertex on left side of domain, (b) outflow vertex on right side of domain}
\label{fig:cornerbc}
\end{center}
\end{figure}

The natural way to specify boundary conditions in a DG scheme is through the boundary fluxes. We have to specify both the fluxes across faces $\hE$ and the vertex fluxes $\tE$. The state outside the domain may be considered as a ghost state and is filled with given boundary condition $\B^*$ so that the same numerical flux as used for interior points can be used on the boundary. At an inflow boundary where $\vel \cdot \un < 0$, the flux $\hE$ is determined from the specified boundary value of $\B^*(x,y,t)$ while at an outflow boundary, it is determined from the interior solution. This is just the upwind principle dictated by the characteristic curves and the numerical flux $\hE$ automatically takes care of this once the ghost value is filled with given boundary condition $\B^*$. In case of the corner flux $\tE$, let us look at an inflow vertex located at the left side of the domain as shown in figure~(\ref{fig:cornerbc}a). Note that we may not have continuity of the normal components across the faces here, e.g., $B_x^* = B_x^U = B_x^D$ may not be satisfied. We will use the data given in the figure and apply the formula~(\ref{eq:nfbk}) to compute the corner flux at inflow boundary. Our numerical experiments show that this leads to a stable scheme and the errors converge at optimal rates even with non-trivial boundary data. At an outflow boundary located on the right side of the domain as shown in figure~(\ref{fig:cornerbc}b), the two outer states are taken to be same as the interior states, $\B_{UR} = \B_{UL}$ and $\B_{DR} = \B_{DL}$, and then the formula formula~(\ref{eq:nfbk}) is used to compute the vertex flux. This is equivalent to computing the flux from the interior state, and can also be written as
\[
\tE = \begin{cases}
E_{DL} & \textrm{if } v_y > 0 \\
E_{UL} & \textrm{otherwise}
\end{cases}
\]
which is the upwind principle based on the characteristics.
%-----------------------------------------------------------------------------------------------
\section{Implementation details}
\label{sec:res}
\begin{figure}
\begin{center}
\includegraphics[width=0.8\textwidth]{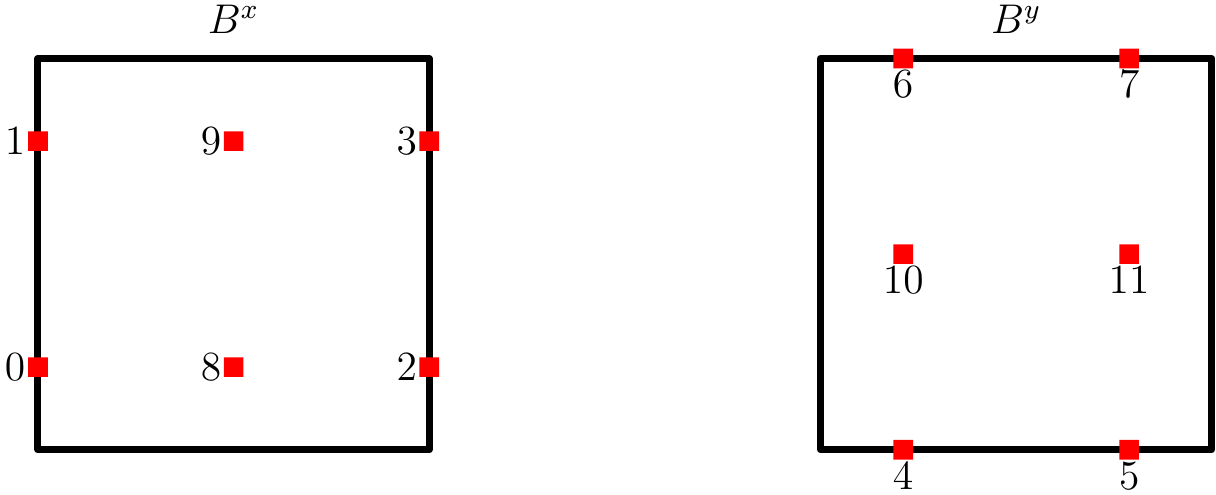}
\caption{Numbering of dofs for $k=1$}
\label{fig:dofs1num}
\end{center}
\end{figure}
The magnetic field $\Bh$ will be approximated in $\rt_k$ in terms of nodal Lagrange polynomials as shown in equation~(\ref{eq:Blag}). As we have seen, the test functions which define the moments are not the same as the trial or basis functions. The test functions for the moments are defined in terms of modal polynomials, each of which has zero mean value except the first one. To define the test functions, we map a cell to the reference cell $[-\half,+\half] \times [-\half,+\half]$ by
\[
\xi = \frac{x - x_0}{\Delta x}, \qquad \eta = \frac{y - y_0}{\Delta y}
\]
where $(x_0,y_0)$ is the cell center and $\Delta x$, $\Delta y$ are the lengths of the sides in the $x$ and $y$ directions, respectively. For degree $k=0$, the test function spaces needed to specify the face moments are
\begin{equation}
P_0(x) = \vspan\{1\}, \qquad P_0(y) = \vspan\{1\}
\label{eq:test0}
\end{equation}
and there are no cell moments in this case. For degree $k=1$, the test function spaces are given by
\begin{equation}
\begin{aligned}
P_1(x) = \vspan\{1, \xi\}, \qquad & P_1(y) = \vspan\{1, \eta\} \\
\partial_x Q_{1,1}(x,y) = \vspan\{1, \eta\}, \qquad & \partial_y Q_{1,1}(x,y) = \vspan\{1, \xi\}
\end{aligned}
\label{eq:test1}
\end{equation}
while for degree $k=2$, they are given by
\begin{equation}
\begin{aligned}
P_2(x) = \vspan\{1, \xi, \xi^2-\otw\}, \qquad P_2(y) = \vspan\{1, \eta, \eta^2-\otw\} \\
\partial_x Q_{2,2}(x,y) = \vspan\{1, \xi, \eta, \xi\eta, \eta^2-\otw, \xi(\eta^2-\otw)\} \\
\partial_y Q_{2,2}(x,y) = \vspan\{1, \xi, \eta, \xi\eta, \xi^2-\otw, (\xi^2-\otw) \eta\}
\end{aligned}
\label{eq:test2}
\end{equation}
It is also possible to use Lagrange polynomials as test functions but we use the above modal test functions in all our computations. We will enumerate the nodal degrees of freedom associated with $\Bxh$, $\Byh$ with a single index. In each cell, there is a local numbering of the dofs. The dofs on the faces are enumerated first in the order $\exl$, $\exr$, $\eyb$, $\eyt$ and then the interior dofs are enumerated. E.g, the case of $k=1$ is illustrated in figure~(\ref{fig:dofs1num}). The set of moment equations~(\ref{eq:mom1})-(\ref{eq:mom4}) leads to a matrix problem where the {\em mass matrix} on each cell has the following structure
\begin{equation}
\begin{bmatrix}
M^x & 0 & 0 & 0 & 0 & 0 \\
0 & M^x & 0 & 0 & 0 & 0 \\
0 & 0 & M^y & 0 & 0 & 0 \\
0 & 0 & 0 & M^y & 0 & 0 \\
N^x_l & N^x_r & 0 & 0 & Q^x & 0 \\
0 & 0 & N^y_b & N^y_t & 0 & Q^y
\end{bmatrix}
\end{equation}
where $M^x$, $M^y$ are $(k+1) \times (k+1)$ matrices arise from the face moments and the remaining matrices arise from the cell moments. The face values are decoupled so that we can solve for the nodal values on each face independently of the other values. Once all the face values are obtained, the interior nodal values can be computed solving the last set of equations which requires inverting the matrices $Q^x$ and $Q^y$. The fifth and sixth rows which correspond to the interior dofs of $B_x$ and $B_y$ are decoupled from one another. Note that since the test functions are different from the basis functions, the mass matrix is somewhat non-standard. E.g., the entries of the first $M^x$ are of the form
\[
\int_{face} (\textrm{basis of $B_x$ with support on left face}) (\textrm{basis of } P_k(y)) \ud y
\]
while the entries of $Q^x$ are of the form
\[
\int_{cell} (\textrm{basis of $B_x$ with interior support}) (\textrm{basis of } \partial_x Q_{k,k}) \ud x \ud y
\]
and the entries of $N^x_l$ are of the form
\[
\int_{cell} (\textrm{basis of $B_x$ with support on left face}) (\textrm{basis of } \partial_x Q_{k,k}) \ud x \ud y
\]
We can write the semi-discrete equations for the $B_x$ components in any cell as
\[
\dd{}{t}[B_x]_l = (M^x)^{-1} R^x_l, \quad \dd{}{t}[B_x]_r = (M^x)^{-1} R^x_r
\]
\[
\dd{}{t}[B_x]_{int} = (Q^x)^{-1} \left[ R^x_{int} - N^x_l \dd{}{t}[B_x]_l - N^x_r \dd{}{t}[B_x]_r \right]
\]
with similar equations for the $B_y$ components. Here $[B_x]_l$, $[B_x]_r$ denote the dofs located on the left and right faces of the cell and $[B_x]_{int}$ denote the interior dofs. We first loop over the all the faces and compute the right hand sides (rhs) of the face moment equations. Then we loop over all the cells and compute the rhs of the cell moment equations. Finally, we can perform one step of the Runge-Kutta scheme.

To give a more concrete view of the scheme and to help the reader to check their own implementation, we give more details about the nodal basis functions for the case of $k=1$, and one can refer to figure~(\ref{fig:dofs1num}) for the following discussion. Following the implementation in \texttt{deal.II}, the nodal basis is defined in terms of the reference cell $[0,1] \times [0,1]$. The nodes on the faces are based on Gauss-Legendre points and are located at $\hxi_0=\half(1-1/\sqrt{3})$ and $\hxi_1=\half(1+1/\sqrt{3})$ on the reference cell. The interior nodes for $\Bxh$ are located at $(\half,\hxi_0)$ and $(\half,\hxi_1)$, while for $\Byh$ are located at $(\hxi_0,\half)$ and $(\hxi_1,\half)$. Define the 1-D Lagrange polynomials
\[
\phi_0(\xi) = \frac{(\xi-\xi_1)(\xi-\xi_2)}{(\xi_0-\xi_1)(\xi_0-\xi_2)}, \quad \phi_1(\xi) = \frac{(\xi-\xi_0)(\xi-\xi_2)}{(\xi_1-\xi_0)(\xi_1-\xi_2)}, \quad \phi_2(\xi) = \frac{(\xi-\xi_0)(\xi-\xi_1)}{(\xi_2-\xi_0)(\xi_2-\xi_1)}
\]
\[
\hphi_0(\xi) = \frac{\xi - \hxi_1}{\hxi_0 - \hxi_1}, \qquad \hphi_1(\xi) = \frac{\xi-\hxi_0}{\hxi_1-\hxi_0}
\]
and the solution can be written as
\begin{eqnarray*}
\Bxh &=& (B_x)_0 \phi_0(\xi) \hphi_0(\eta) + (B_x)_1 \phi_0(\xi) \hphi_1(\eta) + (B_x)_2 \phi_2(\xi) \hphi_0(\eta) + (B_x)_3 \phi_2(\xi) \hphi_1(\eta) + \\
&& (B_x)_8 \phi_1(\xi) \hphi_0(\eta) + (B_x)_9 \phi_1(\xi) \hphi_1(\eta) \\
\Byh &=& (B_y)_4 \hphi_0(\xi) \phi_0(\eta) + (B_y)_5 \hphi_1(\xi) \phi_0(\eta) + (B_y)_6 \hphi_0(\xi) \phi_2(\eta) + (B_y)_7 \hphi_1(\xi) \phi_2(\eta) + \\
&& (B_y)_{10} \hphi_0(\xi) \phi_1(\eta) + (B_y)_{11} \hphi_1(\xi) \phi_1(\eta)
\end{eqnarray*}
where $(B_x)_j$, $(B_y)_j$ are the values at the nodes as numbered in figure~(\ref{fig:dofs1num}). Using the test functions given in~(\ref{eq:test1}), the mass matrix on the reference cell is shown in table~(\ref{tab:mass1}).
\begin{table}
\begin{center}
\[
\left[
\begin{array}{cccccccccccc}
\half & \half & 0 & 0 & 0 & 0 & 0 & 0 & 0 & 0 & 0 & 0 \\
\\
-\frac{1}{4\sqrt{3}} & \frac{1}{4\sqrt{3}} & 0 & 0 & 0 & 0 & 0 & 0 & 0 & 0 & 0 & 0 \\
\\
0 & 0 & \half & \half & 0 & 0 & 0 & 0 & 0 & 0 & 0 & 0 \\
\\
0 & 0 & -\frac{1}{4\sqrt{3}} & \frac{1}{4\sqrt{3}} & 0 & 0 & 0 & 0 & 0 & 0 & 0 & 0 \\
\\
0 & 0 & 0 & 0 & \half & \half & 0 & 0 & 0 & 0 & 0 & 0 \\
\\
0 & 0 & 0 & 0 & -\frac{1}{4\sqrt{3}} & \frac{1}{4\sqrt{3}} & 0 & 0 & 0 & 0 & 0 & 0 \\
\\
0 & 0 & 0 & 0 & 0 & 0 & \half & \half & 0 & 0 & 0 & 0 \\
\\
0 & 0 & 0 & 0 & 0 & 0 & -\frac{1}{4\sqrt{3}} & \frac{1}{4\sqrt{3}} & 0 & 0 & 0 & 0 \\
\\
\otw & \otw & \otw & \otw & 0 & 0 & 0 & 0 & \frac{1}{3} & \frac{1}{3} & 0 & 0 \\
\\
-\frac{\sqrt{3}}{72} & \frac{\sqrt{3}}{72} & -\frac{\sqrt{3}}{72} & \frac{\sqrt{3}}{72} & 0 & 0 & 0 & 0 & -\frac{\sqrt{3}}{72} & \frac{\sqrt{3}}{72} & 0 & 0 \\
\\
0 & 0 & 0 & 0 & \otw & \otw & \otw & \otw & 0 & 0 & \frac{1}{3} & \frac{1}{3} \\
\\
0 & 0 & 0 & 0 & -\frac{\sqrt{3}}{72} & \frac{\sqrt{3}}{72} & -\frac{\sqrt{3}}{72} & \frac{\sqrt{3}}{72} & 0 & 0 & -\frac{\sqrt{3}}{72} & \frac{\sqrt{3}}{72}
\end{array} \right]
\]
\end{center}
\caption{Mass matrix for $k=1$}
\label{tab:mass1}
\end{table}
We perform the integration on the reference cell $[0,1] \times [0,1]$ but the test functions were defined on $[-\half,+\half]\times[-\half,+\half]$, so the coordinates must be transformed as $\xi \to \xi - \half$ and $\eta \to \eta - \half$ before evaluating the test functions. For a general cell, the block matrix for the faces $\expm$ must be scaled by $\Delta y$, those for the faces $\eypm$ must be scaled by $\Delta x$ and the blocks corresponding to the cell moments must be scaled by $\Delta x \Delta y$.

The code is written using \texttt{deal.II}~\cite{BangerthHartmannKanschat2007} which is a \texttt{C++} library that provides building blocks to write finite element programs. In general, the velocity $\vel$ is a function of space and time, and the integrals in the face and cell moment equations have to be computed using some quadrature rule. The integrals on the faces are computed using $(k+2)$-point Gauss-Legendre quadrature and the cell integrals are computed using $(k+2)\times(k+2)$-point Gauss-Legendre quadrature. The time integration is performed by the third order strong stability preserving RK scheme~\cite{Shu1988439}. The time step is chosen according to the following condition
\[
\Delta t = \frac{\cfl}{(2k+1)\max\left(\frac{|v_x|}{\Delta x} + \frac{|v_y|}{\Delta y}\right)}
\]
and in all the test cases, we choose $\cfl = 0.8$. The above formula is motivated by the time step restrictions normally used in DG schemes and in all our tests, we have found that the above choice was stable. Of course, a more rigorous stability analysis has to be performed in a future work.
%-----------------------------------------------------------------------------------------------
\section{Numerical results}
In this section, we provide numerical evidence to the approximation of vector fields using Raviart-Thomas polynomials. We then use these polynomials to solve induction equation for cases with zero and non-zero divergence, and numerically show that optimal convergence orders are achieved.
%-----------------------------------------------------------------------------------------------
\subsection{Test 1: Approximation of smooth fields}
In this section, we test the accuracy of projecting a given divergence-free field onto the polynomial space $\rt_k$ using the moments. For the first example titled Test 1a, we use a divergence-free field given by $\B = (\partial_y \Phi, -\partial_x \Phi)$ where
\[
\Phi(x,y) = \sin(2\pi x) \sin(2\pi y), \qquad (x,y) \in [0,1] \times [0,1]
\]
For the second example titled Test 1b, we take a divergent field given by $\B = \nabla \Phi$ where
\[
\Phi(x,y) = \frac{1}{10} \exp[-20(x^2 + y^2)], \qquad (x,y) \in [-1,+1] \times [-1,+1]
\]
The error for the first example are shown in tables~(\ref{tab:app11}),~(\ref{tab:app12}) which shows the optimal convergence rates consistent with the error estimates given in equation~(\ref{eq:app1}). The norm of the divergence is small; the variation seen is due to the difficulty in accurately computing a quantity that is zero due to roundoff errors. The errors for second example are shown in tables~(\ref{tab:app21}),~(\ref{tab:app22}) and we observe that both the function and its divergence converge at the same optimal rate.
\begin{table}
\begin{center}
\begin{tabular}{|c|c|c|c|} \hline
$h$ & 
\multicolumn{2}{|c|}{$\|B-B_h\|_{L^2(\Omega)}$} & $\|div(B_h)\|_{L^2(\Omega)}$\\ \hline
0.1250 & 1.0189e-01 & - & 3.7147e-14\\ \hline
0.0625 & 2.5519e-02 & 2.00 & 9.5162e-14\\ \hline
0.0312 & 6.3826e-03 & 2.00 & 3.7880e-13\\ \hline
0.0156 & 1.5958e-03 & 2.00 & 1.4840e-12\\ \hline
0.0078 & 3.9896e-04 & 2.00 & 5.8016e-12\\ \hline
\end{tabular}
\caption{Test 1a: Approximation error convergence for $k=1$}
\label{tab:app11}
\end{center}
\end{table}

\begin{table}
\begin{center}
\begin{tabular}{|c|c|c|c|} \hline
$h$ & 
\multicolumn{2}{|c|}{$\|B-B_h\|_{L^2(\Omega)}$} & $\|div(B_h)\|_{L^2(\Omega)}$\\ \hline
0.1250 & 6.7521e-03 & - & 1.3265e-13\\ \hline
0.0625 & 8.4659e-04 & 3.00 & 3.7389e-13\\ \hline
0.0312 & 1.0590e-04 & 3.00 & 1.3266e-12\\ \hline
0.0156 & 1.3241e-05 & 3.00 & 5.2716e-12\\ \hline
0.0078 & 1.6552e-06 & 3.00 & 2.0924e-11\\ \hline
\end{tabular}
\caption{Test 1a: Approximation error convergence for $k=2$}
\label{tab:app12}
\end{center}
\end{table}

\begin{table}
\begin{center}
\begin{tabular}{|c|c|c|c|c|} \hline
$h$ & 
\multicolumn{2}{|c|}{$\|B-B_h\|_{L^2(\Omega)}$} & 
\multicolumn{2}{|c|}{$\|div(B)-div(B_h)\|_{L^2(\Omega)}$}\\ \hline
0.0625 & 9.0930e-04 & - & 2.7438e-02 & -\\ \hline
0.0312 & 2.2445e-04 & 2.02 & 6.9076e-03 & 1.99\\ \hline
0.0156 & 5.5927e-05 & 2.00 & 1.7299e-03 & 2.00\\ \hline
0.0078 & 1.3970e-05 & 2.00 & 4.3267e-04 & 2.00\\ \hline
0.0039 & 3.4918e-06 & 2.00 & 1.0818e-04 & 2.00\\ \hline
\end{tabular}
\caption{Test 1b: Approximation error convergence for $k=1$}
\label{tab:app21}
\end{center}
\end{table}

\begin{table}
\begin{center}
\begin{tabular}{|c|c|c|c|c|} \hline
$h$ & 
\multicolumn{2}{|c|}{$\|B-B_h\|_{L^2(\Omega)}$} & 
\multicolumn{2}{|c|}{$\|div(B)-div(B_h)\|_{L^2(\Omega)}$}\\ \hline
0.0625 & 4.7750e-05 & - & 1.8703e-03 & -\\ \hline
0.0312 & 5.9190e-06 & 3.01 & 2.3550e-04 & 2.99\\ \hline
0.0156 & 7.3827e-07 & 3.00 & 2.9491e-05 & 3.00\\ \hline
0.0078 & 9.2233e-08 & 3.00 & 3.6881e-06 & 3.00\\ \hline
0.0039 & 1.1528e-08 & 3.00 & 4.6106e-07 & 3.00\\ \hline
\end{tabular}
\caption{Test 1b: Approximation error convergence for $k=2$}
\label{tab:app22}
\end{center}
\end{table}
%-----------------------------------------------------------------------------------------------
\subsection{Test 2: Smooth test case, divergence-free solution}
The initial condition is given by $\B_0 = (\partial_y \Phi, -\partial_x \Phi)$ where
\[
\Phi(x,y) = \frac{1}{10} \exp[-20((x-1/2)^2 + y^2)]
\]
and the velocity field is $\vel = (y, -x)$. The exact solution is a pure rotation of the initial condition and is given by
\[
\B(\vr,t) = R(t) \B_0(R(-t) \vr), \qquad R(t) = \begin{bmatrix}
\cos t & -\sin t \\
\sin t & \cos t \end{bmatrix}
\]
By construction, the exact solution has zero divergence initially and hence at future times also. We solve this problem on two domain sizes which helps to show that the method is able to implement non-trivial boundary conditions that depend on both space and time in a stable and accuracy preserving manner.
\subsubsection{Test 2a: Large domain}
\begin{figure}
\begin{center}
\begin{tabular}{ccc}
\includegraphics[width=0.32\textwidth]{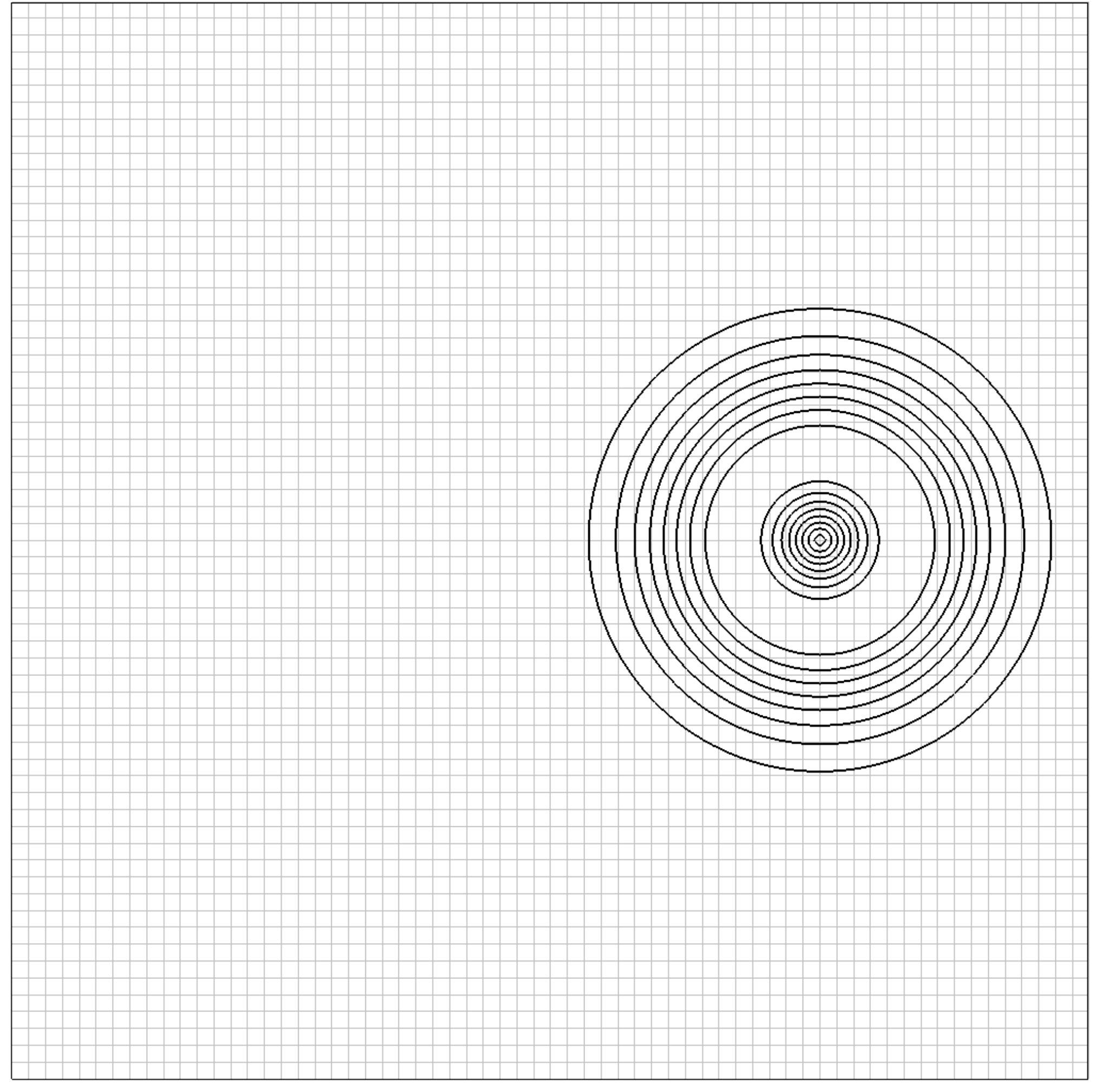} &
\includegraphics[width=0.32\textwidth]{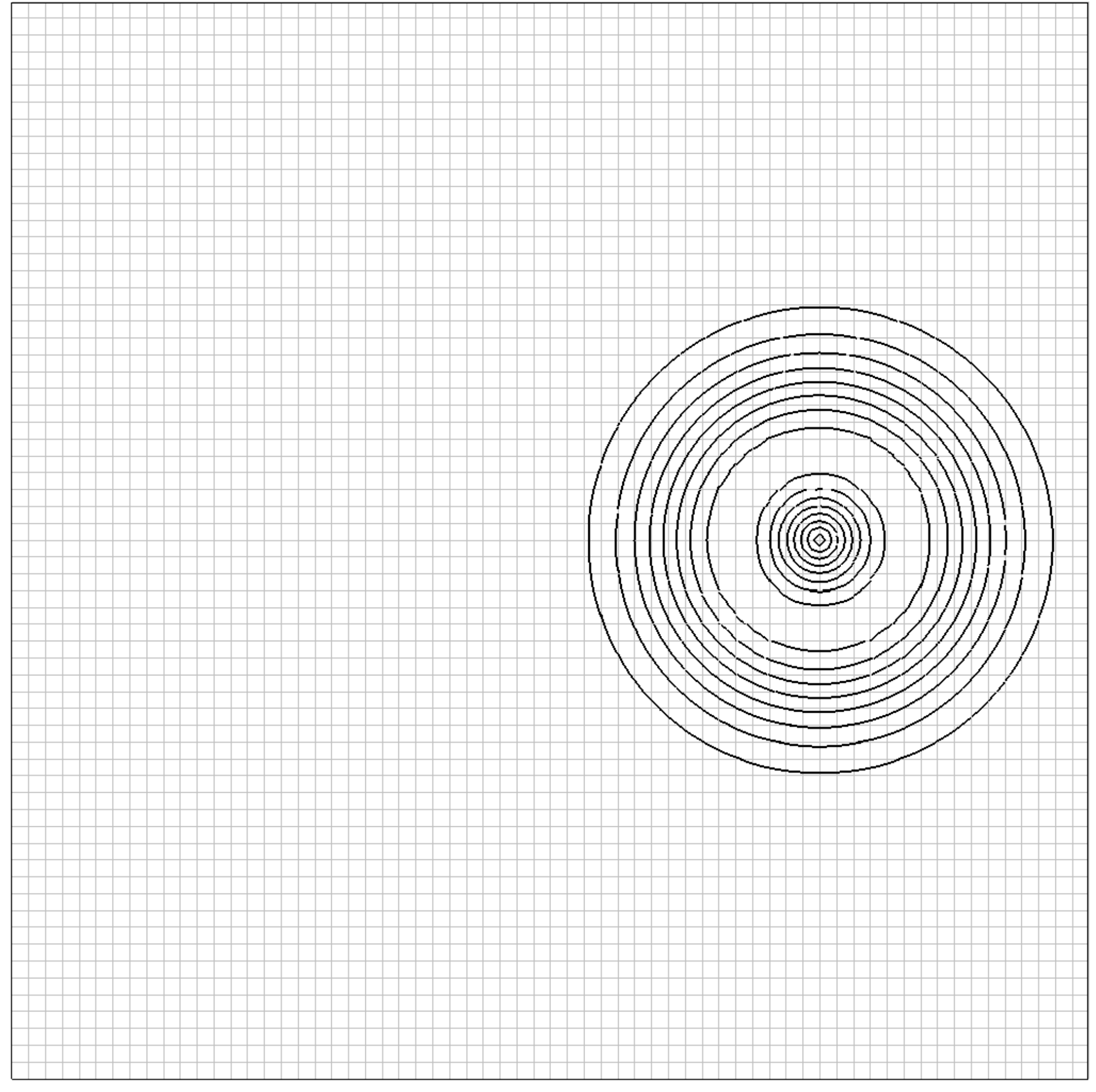} &
\includegraphics[width=0.32\textwidth]{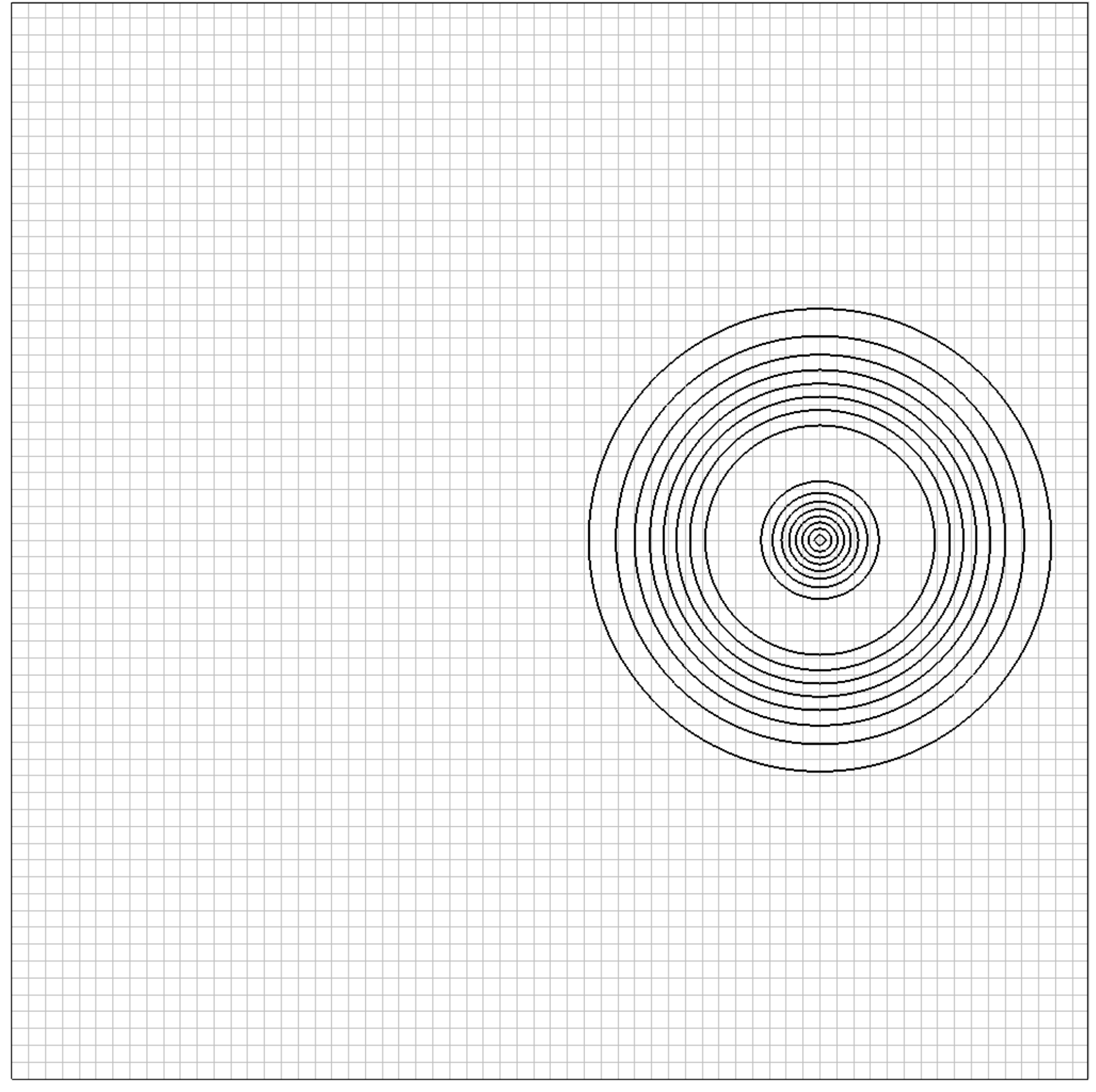} \\
(a) & (b) & (c)
\end{tabular}
\end{center}
\caption{Contour of $|\B_h|$ for Test 2a, 10 contours between 0 and 0.3867: (a) initial, (b) final, $k=1$, (c) final, $k=2$}
\label{fig:anal2}
\end{figure}

We compute the numerical solution on the computational domain $[-1,+1] \times [-1,+1]$ upto a final time of $T=2\pi$ at which time the solution comes back to the initial condition. At the boundary, the solution is nearly zero due to exponential decay of the solution. Figure~(\ref{fig:anal2}) shows the contours of the solution at the final time and the mesh of $64 \times 64$ cells used in this simulation is inlaid in the background. We clearly see the improvement in the solution when we go from $k=1$ to $k=2$, corresponding to second and third order schemes respectively. The problem is solved on a sequence of refined meshes and the corresponding error norms are shown in tables~(\ref{tab:anal2ak1}) and (\ref{tab:anal2ak2}), which shows the design order of accuracy is being achieved.

\begin{table}
\begin{center}
\begin{tabular}{|c|c|c|c|}
\hline
$h$ & \multicolumn{2}{|c|}{$\norm{\B_h - \B}_{L^2(\Omega)}$} & $\|div(\B_h)\|_{L^2(\Omega)}$ \\
\hline
0.0312 & 2.1427e-03 & - & 6.0137e-14 \\
0.0156 & 3.2571e-04 & 2.71 & 1.8566e-13 \\
0.0078 & 5.9640e-05 & 2.45 & 5.8486e-13 \\
0.0039 & 1.3209e-05 & 2.17 & 1.8853e-12 \\
\hline
\end{tabular}
\caption{Convergence of error for Test 2a with $k=1$}
\label{tab:anal2ak1}
\end{center}
\end{table}

\begin{table}
\begin{center}
\begin{tabular}{|c|c|c|c|}
\hline
$h$ & \multicolumn{2}{|c|}{$\norm{\B_h - \B}_{L^2(\Omega)}$} & $\|div(\B_h)\|_{L^2(\Omega)}$ \\
\hline
0.0625 & 2.4003e-04 & - & 4.9081e-14 \\
0.0312 & 2.5212e-05 & 3.25 & 1.4299e-13 \\
0.0156 & 3.0946e-06 & 3.02 & 4.5663e-13 \\
0.0078 & 3.8448e-07 & 3.00 & 1.5058e-12 \\
\hline
\end{tabular}
\caption{Convergence of error for Test 2a with $k=2$}
\label{tab:anal2ak2}
\end{center}
\end{table}

\subsubsection{Test 2b: Small domain}
\begin{figure}
\begin{center}
\begin{tabular}{ccc}
\includegraphics[width=0.32\textwidth]{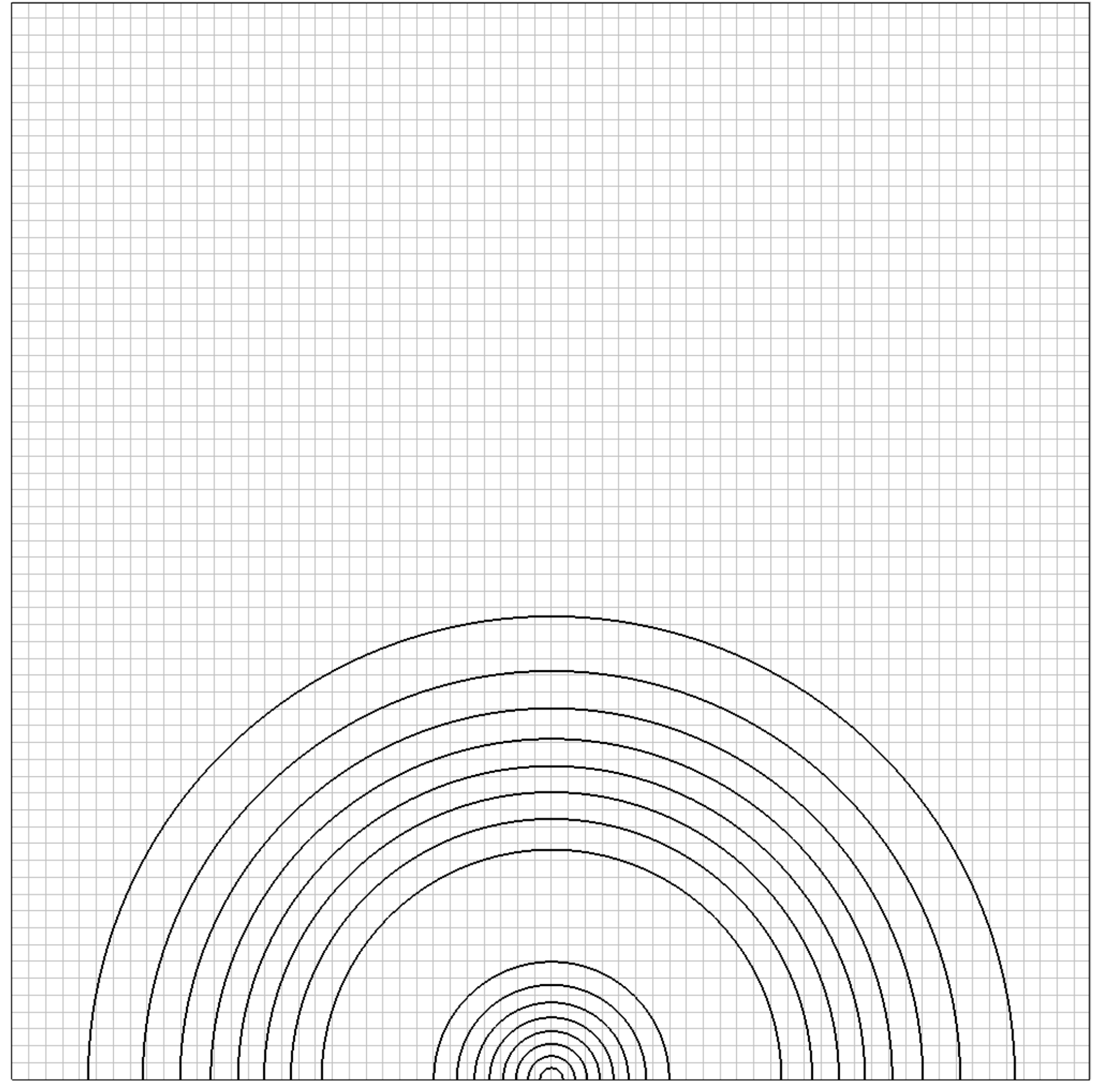} &
\includegraphics[width=0.32\textwidth]{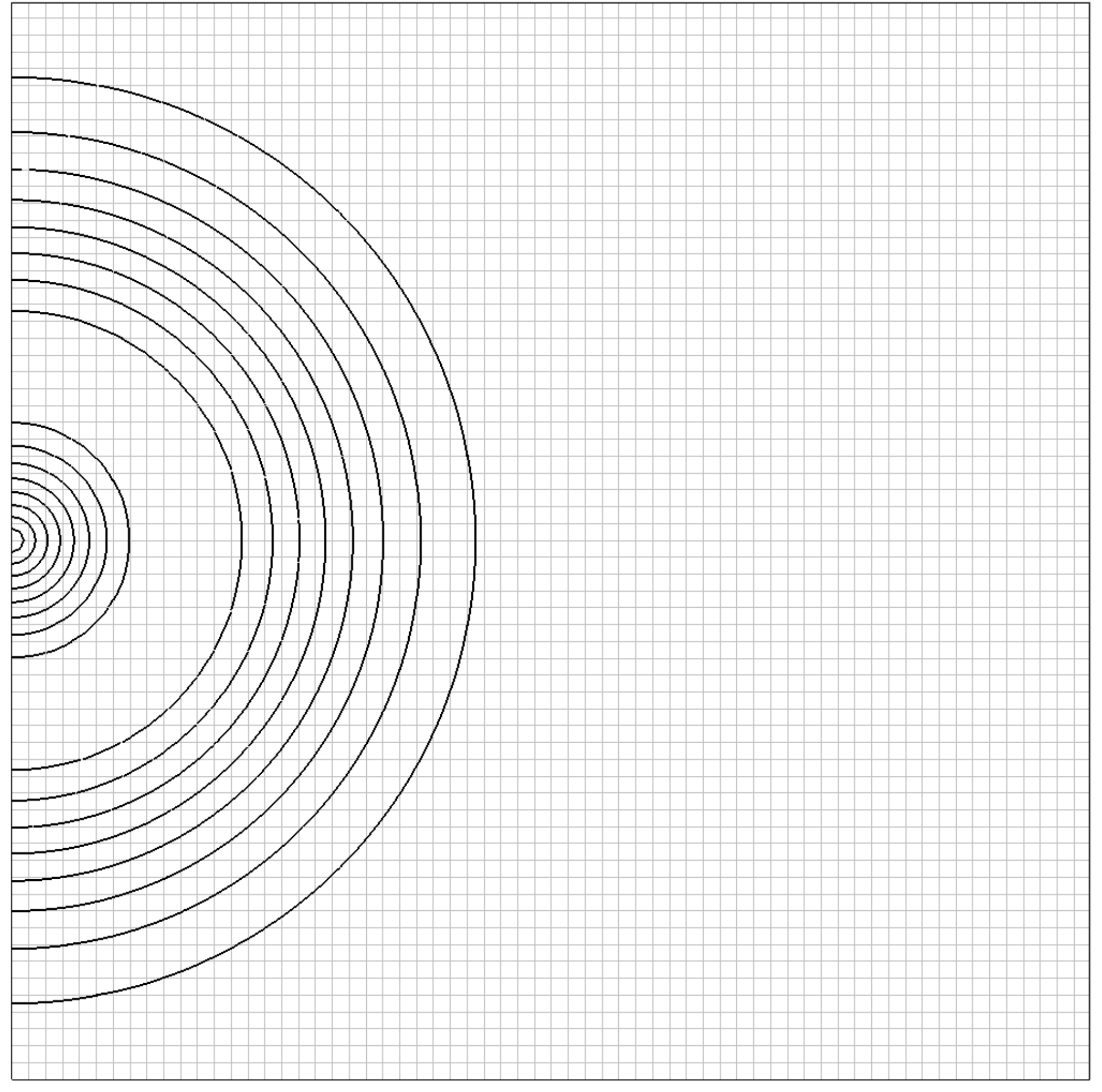} &
\includegraphics[width=0.32\textwidth]{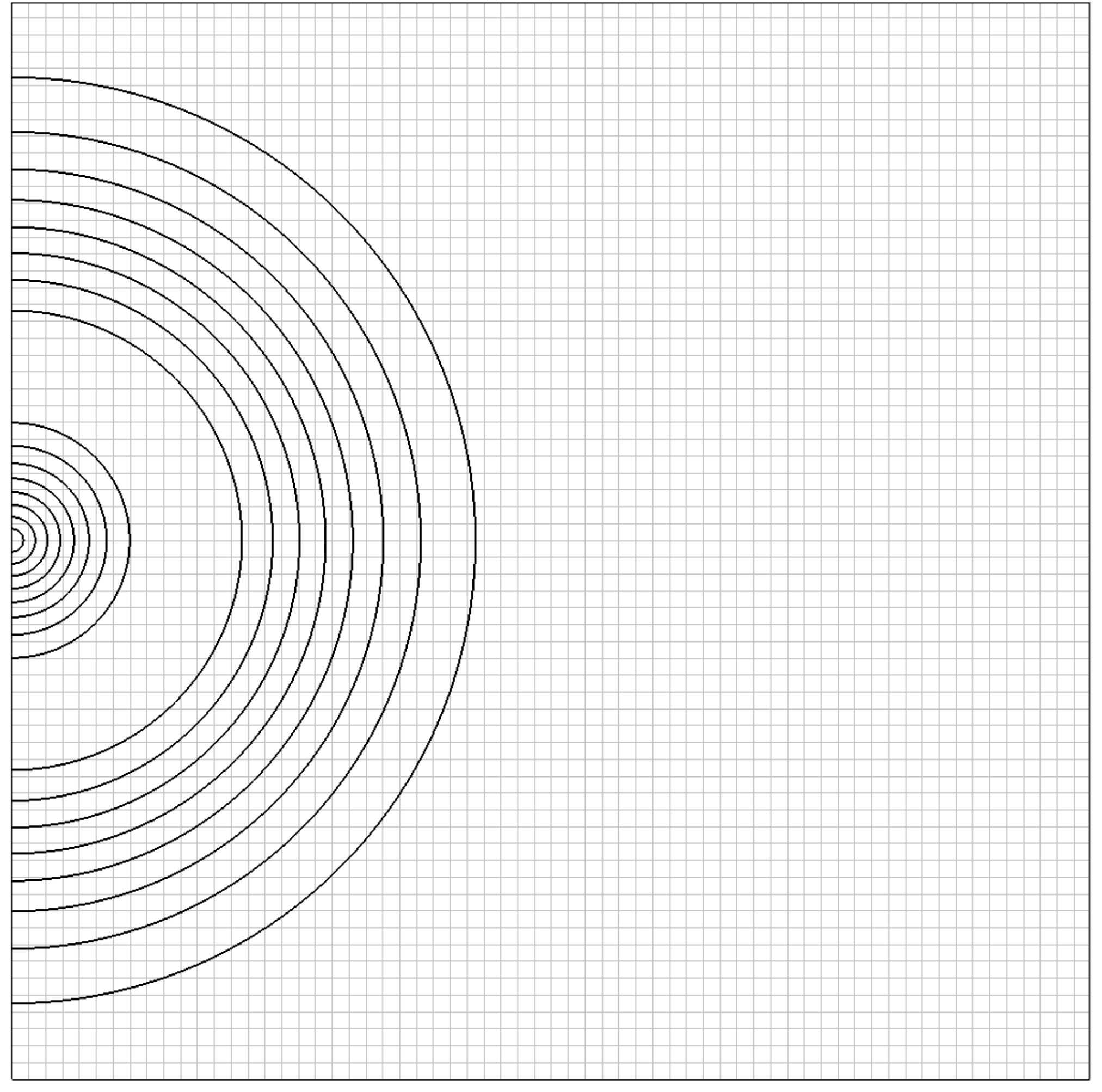} \\
(a) & (b) & (c)
\end{tabular}
\end{center}
\caption{Contour of $|\B_h|$ for Test 2b, 10 contours between 0 and 0.3867: (a) initial, (b) final, $k=1$, (c) final, $k=2$}
\label{fig:anal2b}
\end{figure}

We compute the numerical solution on the computational domain $[0,1] \times [0,1]$ upto a final time of $T=\pi/2$. Due to this finite domain, the solution at the boundary is non-trivial. Figure~(\ref{fig:anal2b}) shows sample solution on a grid of $64 \times 64$ cells; the initial solution profile is located at the lower part of the domain and at the final time, this has rotated by 90 degrees in counter-clockwise direction and part of the solution has exited the domain. We also compute the solution on a sequence of successively refined grids and the error is shown in tables~(\ref{tab:anal2bk1}) and (\ref{tab:anal2bk2}) respectively for the second and third order cases. These results indicate that the design order of accuracy has been achieved even in the presense of non-trivial boundary conditions.

\begin{table}
\begin{center}
\begin{tabular}{|c|c|c|c|}
\hline
$h$ & \multicolumn{2}{|c|}{$\norm{\B_h - \B}_{L^2(\Omega)}$} & $\|div(\B_h)\|_{L^2(\Omega)}$ \\
\hline
0.0312 & 6.5882e-04 & - & 2.8687e-14 \\
0.0156 & 1.4979e-04 & 2.13 & 9.8666e-14 \\
0.0078 & 3.6394e-05 & 2.04 & 3.2902e-13 \\
0.0039 & 9.0308e-06 &  2.01 & 1.1356e-12 \\
\hline
\end{tabular}
\caption{Convergence of error for Test 2b with $k=1$}
\label{tab:anal2bk1}
\end{center}
\end{table}

\begin{table}
\begin{center}
\begin{tabular}{|c|c|c|c|}
\hline
$h$ & \multicolumn{2}{|c|}{$\norm{\B_h - \B}_{L^2(\Omega)}$} & $\|div(\B_h)\|_{L^2(\Omega)}$ \\
\hline
0.0625 & 1.4110e-04 & - & 2.4986e-14 \\
0.0312 & 1.7238e-05 & 3.03 & 7.9129e-14 \\
0.0156 & 2.1442e-06 &  3.00 & 2.5910e-13 \\
0.0078 & 2.6749e-07 & 3.00 & 9.2720e-13 \\
\hline
\end{tabular}
\caption{Convergence of error for Test 2b with $k=2$}
\label{tab:anal2bk2}
\end{center}
\end{table}
%-----------------------------------------------------------------------------------------------
\subsection{Test 3: Smooth test case, divergent solution}
\begin{figure}
\begin{center}
\begin{tabular}{ccc}
\includegraphics[width=0.32\textwidth]{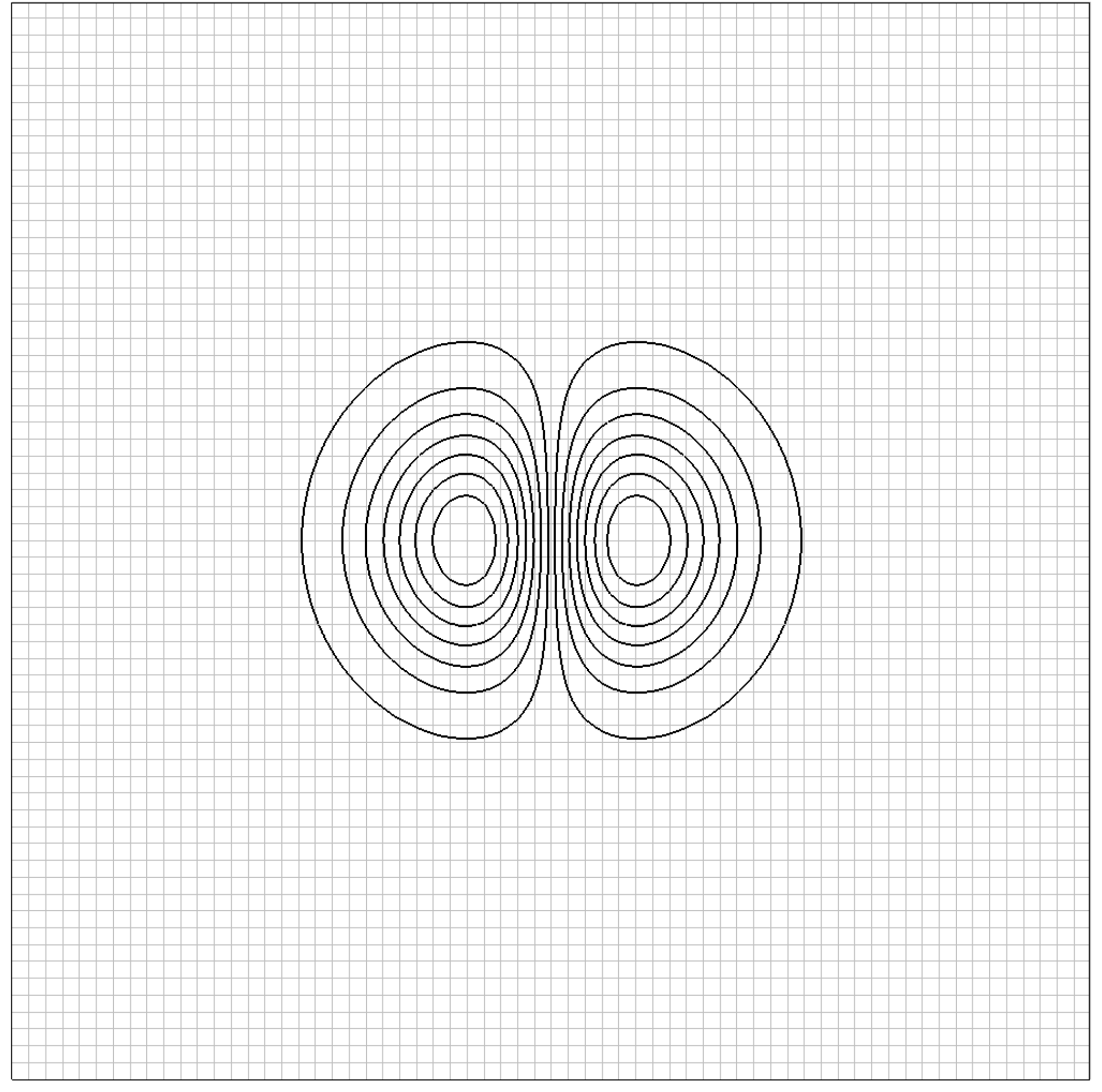} &
\includegraphics[width=0.32\textwidth]{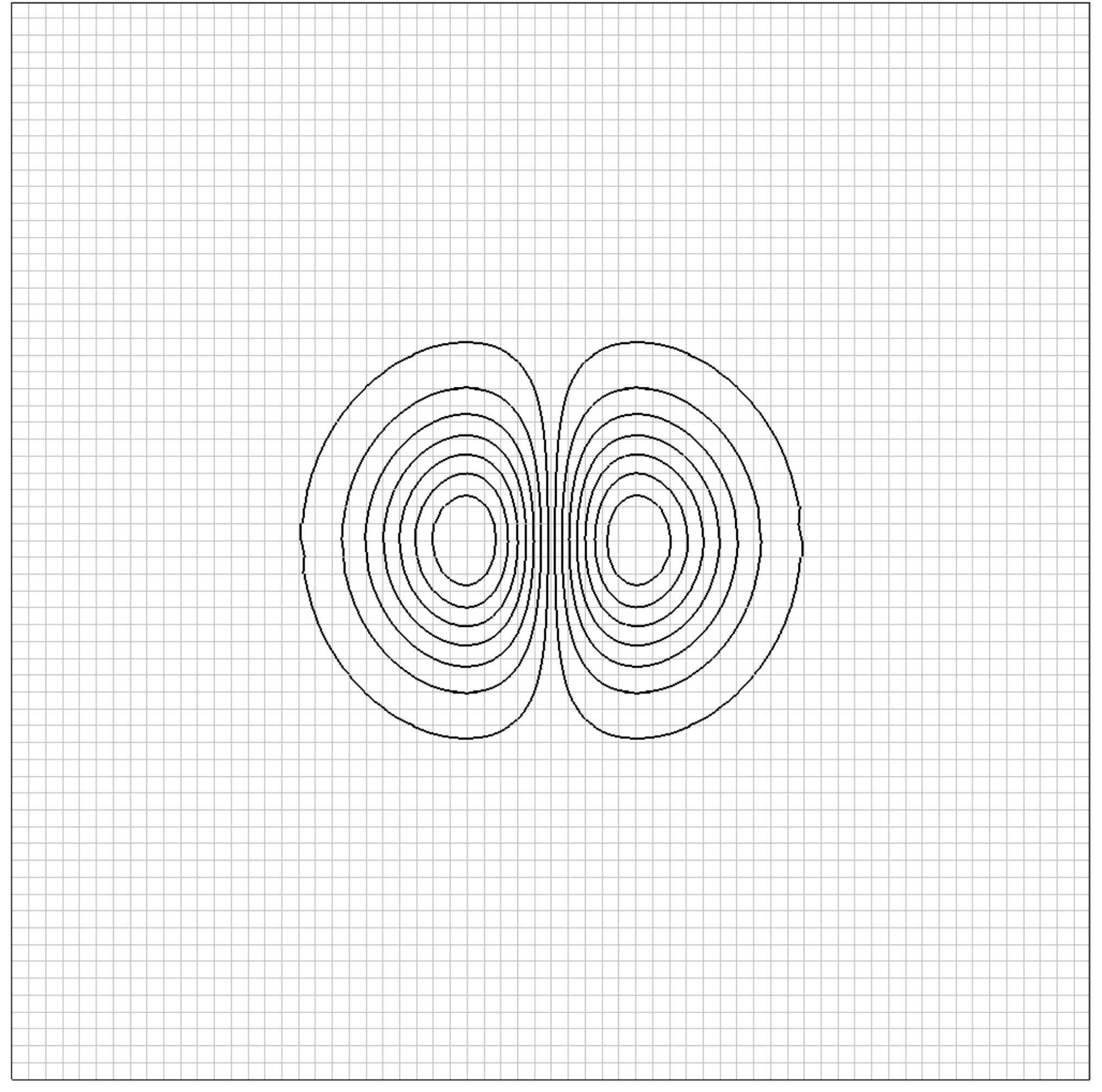} &
\includegraphics[width=0.32\textwidth]{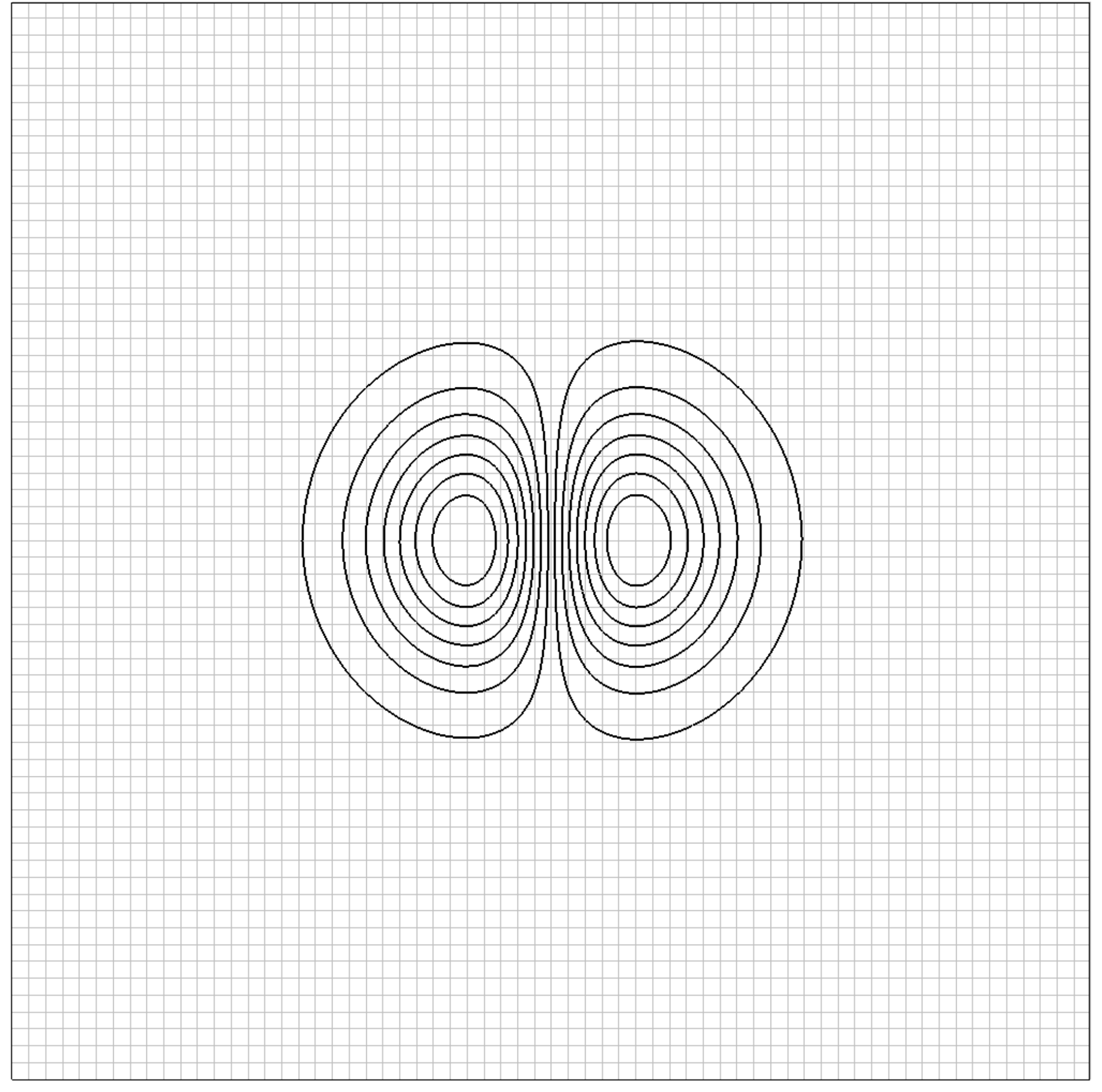} \\
(a) & (b) & (c)
\end{tabular}
\end{center}
\caption{Contour plot of $B_x$ for Test 3 showing 16 contours between -0.3838 and +0.3838: (a) initial condition, (b) final, $k=1$ and (c) final, $k=2$}
\label{fig:anal4}
\end{figure}

In this problem, we generate an exact solution by the method of manufactured solutions. The exact solution is taken to be
\[
\B(x,y,t) = \begin{bmatrix}
\cos t & -\sin t \\
\sin t & \cos t \end{bmatrix} \B_0(x,y)
\]
where
\[
\B_0 = \nabla \phi, \qquad \phi = \frac{1}{10} \exp(-20(x^2 + y^2))
\]
and the velocity field is taken as
\[
\vel = \nabla^\top \psi, \qquad \psi = \frac{1}{\pi} \sin(\pi x) \sin(\pi y)
\]
Note that by construction, the solution has non-zero divergence. The right hand side source term $\M$ is computed from the above solution using the formula $\M = -\df{\B}{t} + \nabla\times(\vel \times \B)$. The problem is solved on the domain $[-1,+1] \times [-1,+1]$ until a final time of $T = 2\pi$. The $x$ component of the solution on a grid of $64 \times 64$ cells is shown in figure~(\ref{fig:anal4}) for the second and third order schemes. The solution contours rotate around the origin and at the final time, the contours should coincide with the initial condition. The figures show very similar contours at the final time and we see the third order being slightly better. The convergence of the error in the solution and its divergence is shown in tables~(\ref{tab:anal4k1}), (\ref{tab:anal4k2}), respectively for the case of $k=1$ and $k=2$. We see that both the solution and its divergence converge at the optimal rate of $k+1$. This shows that in case of CED, we can compute the charge density also to optimal accuracy since it depends on the divergence of the solution.

\begin{table}
\begin{center}
\begin{tabular}{|c|c|c|c|c|c|}
\hline
$h$ & \multicolumn{2}{|c|}{$\norm{\B_h - \B}_{L^2(\Omega)}$} & \multicolumn{2}{|c|}{$\|div(\B) -div(\B_h)\|_{L^2(\Omega)}$} \\
\hline
0.0312 & 8.5550e-04 & - & 6.9076e-03 & - \\
0.0156 & 1.8915e-04 & 2.17 & 1.7299e-03 & 1.99 \\
0.0078 & 3.8730e-05 & 2.29 & 4.3267e-04 & 1.99 \\
0.0039 & 7.8346e-06 & 2.30 & 1.0818e-04 & 1.99 \\
\hline
\end{tabular}
\caption{Convergence of error for Test 3 with $k=1$}
\label{tab:anal4k1}
\end{center}
\end{table}

\begin{table}
\begin{center}
\begin{tabular}{|c|c|c|c|c|c|}
\hline
$h$ & \multicolumn{2}{|c|}{$\norm{\B_h - \B}_{L^2(\Omega)}$} & \multicolumn{2}{|c|}{$\|div(\B) -div(\B_h)\|_{L^2(\Omega)}$} \\
\hline
0.0625 & 3.4775e-04 &-& 1.8703e-03 & - \\
0.0312 & 3.3408e-05 & 3.38 & 2.3550e-04 & 2.99 \\
0.0156 & 3.0287e-06 & 3.46 & 2.9491e-05 & 2.99 \\
0.0078 & 2.7345e-07 & 3.47 & 3.6881e-06 & 2.99 \\
\hline
\end{tabular}
\caption{Convergence of error for Test 3 with $k=2$}
\label{tab:anal4k2}
\end{center}
\end{table}
%-----------------------------------------------------------------------------------------------
\subsection{Test 4: Discontinuous test case}
\begin{figure}
\begin{center}
\begin{tabular}{cc}
\includegraphics[width=0.45\textwidth]{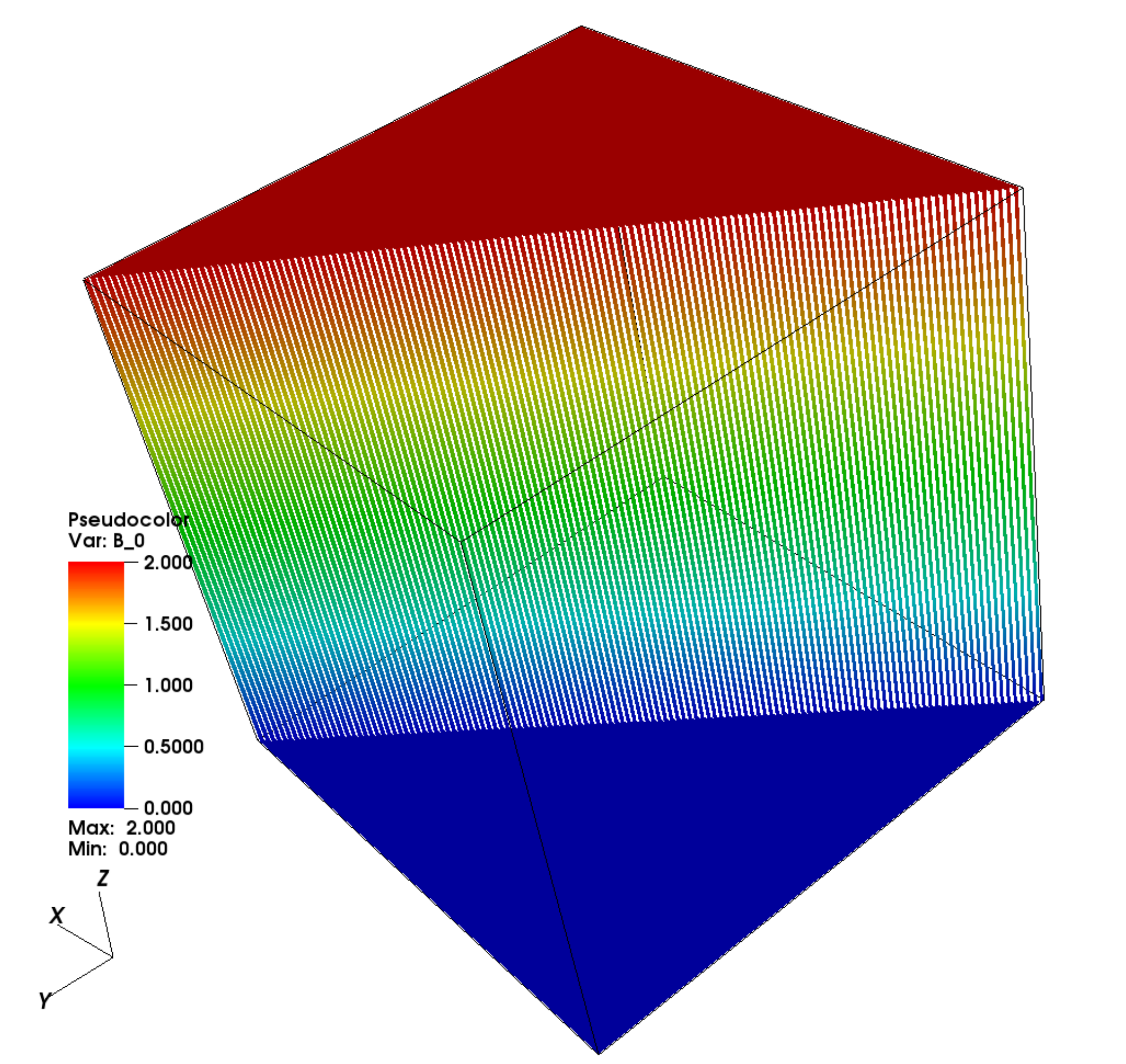} &
\includegraphics[width=0.45\textwidth]{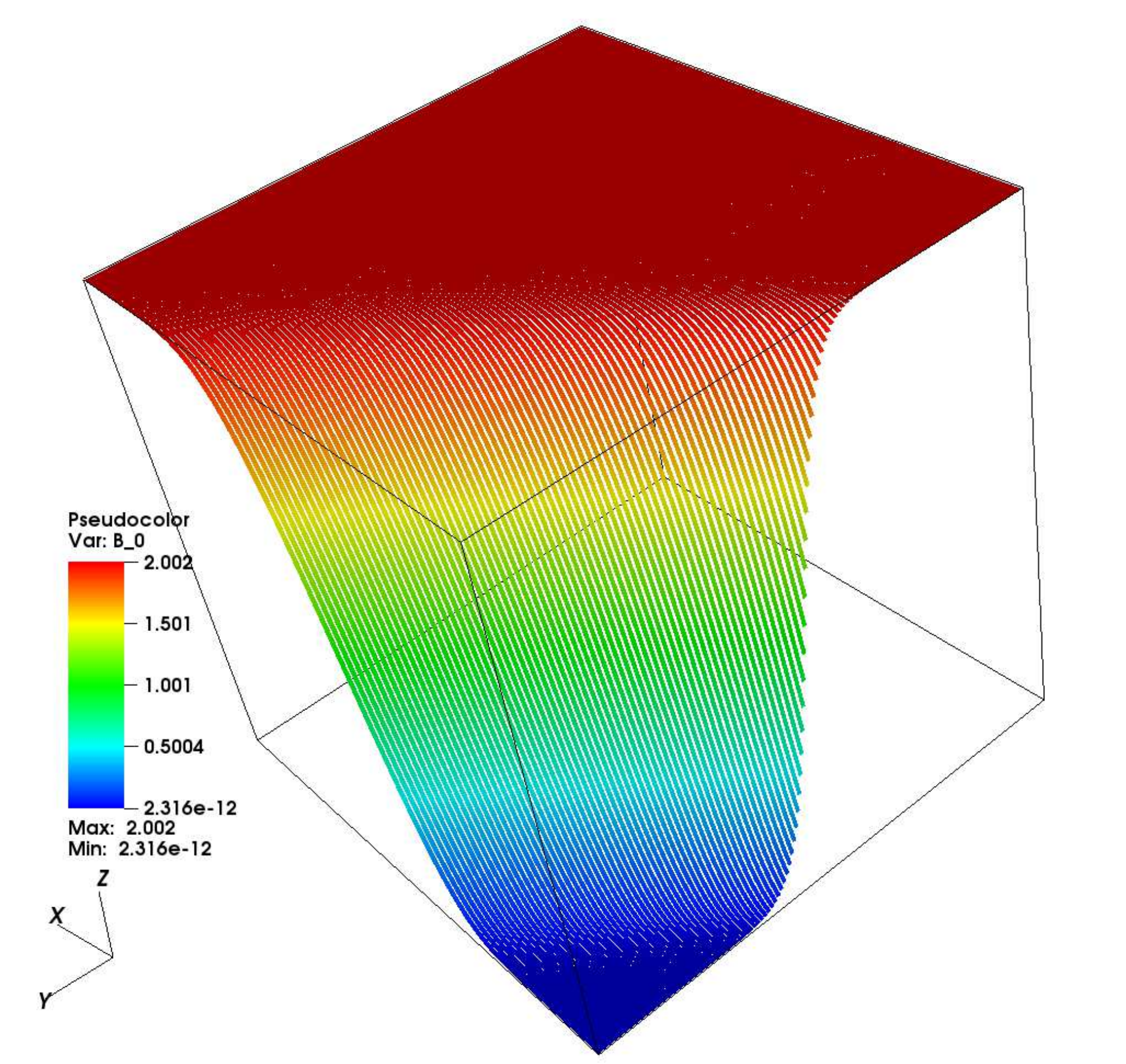} \\
(a) & (b) \\
\includegraphics[width=0.45\textwidth]{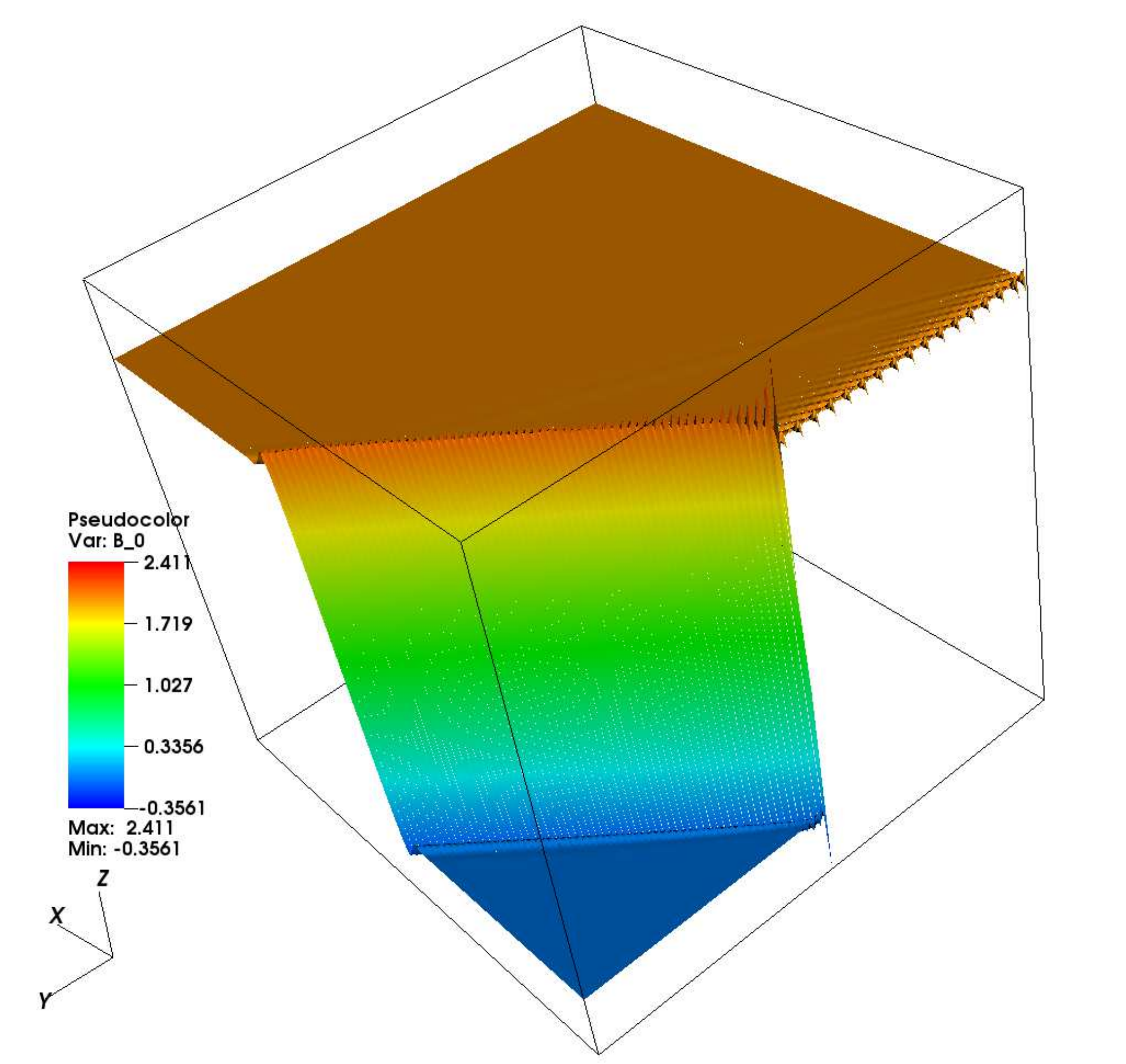} &
\includegraphics[width=0.45\textwidth]{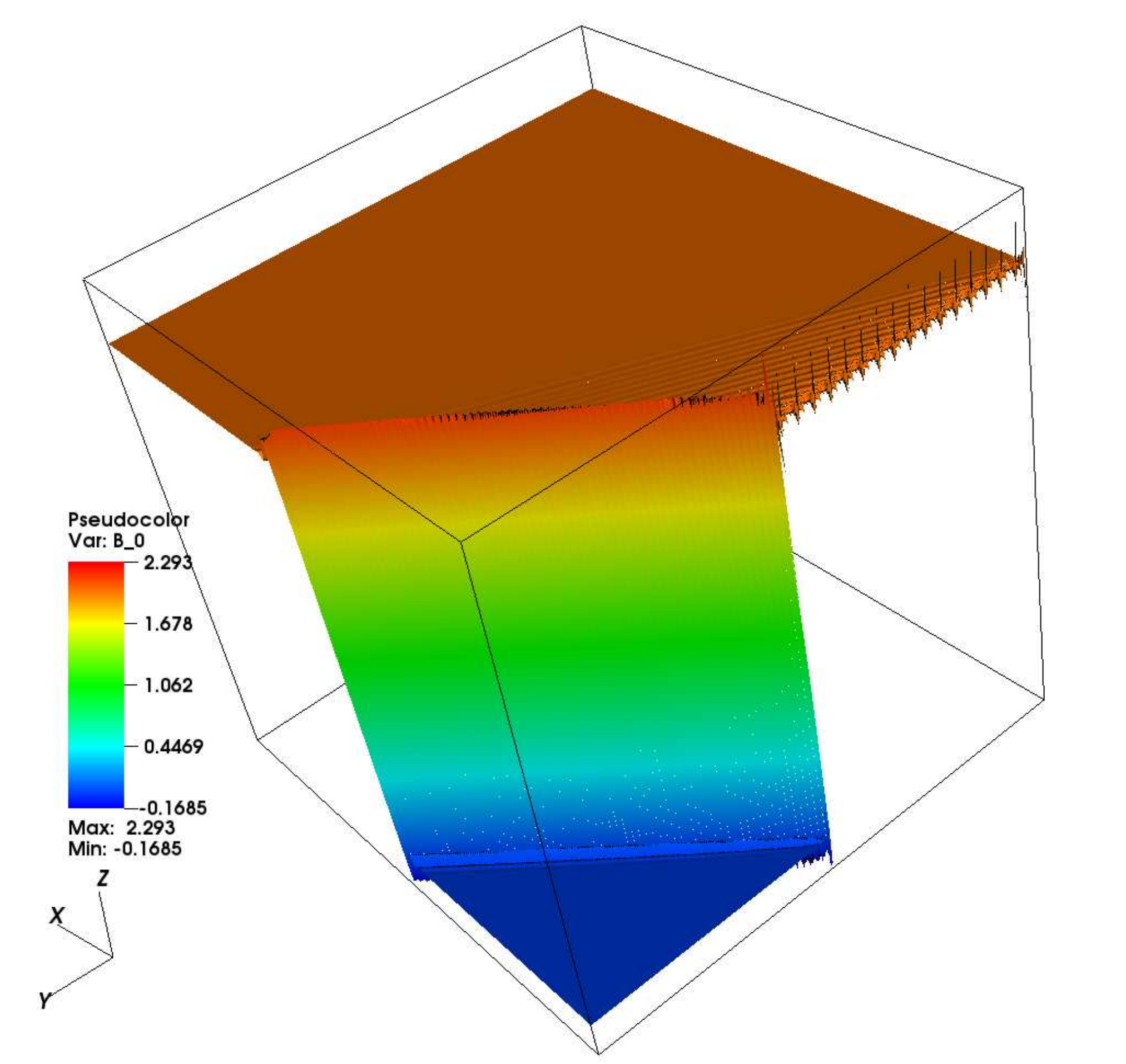} \\
(c) & (d)
\end{tabular}
\caption{Solution for Test 4 on grid of $128 \times 128$ cells: (a) initial condition, (b) $k=0$, (c) $k=1$, (d) $k=2$}
\label{fig:anal3}
\end{center}
\end{figure}

The scheme developed so far is not suitable for computing discontinuous solutions since we need some form of limiting to control the Gibbs oscillations. However, due to the discontinuous Galerkin and upwind nature of the scheme, it should still be stable for a linear PDE like the induction equation in the sense that the computations should not blow up and any oscillations should be restricted to regions close to the discontinuities. We will show in this test case that the scheme indeed achieves these objectives. We take the potential
\[
\Phi(x,y) = \begin{cases}
2y - 2x & \textrm{if } x > y \\
0 & \textrm{otherwise}
\end{cases}
\]
and the velocity field is $\vel = (1,2)$. This leads to a discontinuous magnetic field with the discontinuity along the line $x=y$ and the initial magnetic field is given by
\[
\B_0 = \begin{cases}
(2,2) & \textrm{if } x > y \\
(0,0) & \textrm{if } x < y
\end{cases}
\]
The exact solution is obtained by a translation of the initial condition and is given by
\[
\B(x,y,t) = \B_0(x-t, y-2t)
\]
We compute this solution on a grid of $128 \times 128$ cells upto a final time of $T=0.5$ units. The solutions for degree $k=0$, $k=1$ and $k=2$ are shown figure~(\ref{fig:anal3}) in terms of surface plots of the $x$ component of $\B_h$. For $k=0$, the solution is non-oscillatory and corresponds to a first order scheme. Note that even though this corresponds to linear variation inside the cell, the solutions are non-oscillatory as we expect from a first order method. For higher order schemes, we see that there are oscillations around the discontinuity line and also near the inlet portion where the discontinuity hits the boundary. In other regions we do not observe the spread of these oscillations which indicates the DG scheme has a stabilizing effect. At this final time, the divergence norm of the solution is 3.9055e-13, 2.7616e-12 and 8.1331e-12 respectively for the three cases, showing that even in this case the divergence-free property is maintained by the scheme.

%-----------------------------------------------------------------------------------------------
\section{Summary and Conclusions}
A new type of DG scheme has been proposed to deal with problems involving a divergence constraint by utilizing approximating polynomials spaces based on Raviart-Thomas polynomials. These polynomials naturally provide divergence-free approximations on any mesh provided the data satisfies this property. By carefully evolving the moments used to construct the approximating polynomials by a DG scheme, we are able to preserve the divergence-free property at future times also without any extra reconstruction process or modification of solution. When the divergence is non-zero, it is computed to same accuracy as the solution which is useful to approximate the electric field in Maxwell equations. The DG schemes require multi-dimensional fluxes which have been recently proposed in the literature for various systems like MHD and CED. The use of fluxes obtained from a Riemann solver is seen to lead to stable schemes for induction equation that show optimal convergence rates in numerical tests. The present paper is devoted to the mathematical aspects of satisfying the involution constraint inherent in Faraday's law with the help of a specially-formulated DG scheme. In future work we will show how this synthesis between DG schemes and multidimensional Riemann solvers yields superior DG schemes for several involution-constrained systems like MHD and CED with extensions to unstructured, isoparametric and adaptive grids.
%-----------------------------------------------------------------------------------------------
\section*{Acknowledgments}
The author would like to acknowledge the support received from the Airbus Chair on Mathematics of Complex Systems established at TIFR-CAM by the Airbus Foundation for carrying out this work. The author also thanks Dinshaw S. Balsara for many discussions which were helpful in formulating these ideas. Finally, the author would like to thank the anonymous reviewer whose comments helped to improve the presentation of the paper.
%-----------------------------------------------------------------------------------------------

\end{document}